\pdfoutput=1
\documentclass{article}

\usepackage{makeidx}         
\usepackage{graphicx}        
\usepackage{multicol}        
\usepackage[bottom]{footmisc}
\usepackage{cite} 
 \usepackage{amssymb}
 \usepackage{amsmath}
 \usepackage{latexsym}
\usepackage{url}

   \newtheorem{theorem}{Theorem}
   \newtheorem{proposition}{Proposition}
   \newtheorem{corollary}{Corollary}

  \newtheorem{definition}{Definition}
  \numberwithin{equation}{section}

\newcommand{\qed}{\hspace*{\fill} $\Box$ \ifmmode \else
    \par\addvspace\topsep\fi}
\newenvironment {proof}{\par\addvspace\topsep\noindent {\it Proof.}
    \ignorespaces }

\makeindex             

\title{Graph Polynomials and Their Applications I: \\The Tutte Polynomial} 
\author{Joanna A. Ellis-Monaghan\,$^1$
        \and
        Criel Merino\,$^2$
        }
\date{}

\begin{document}
 \maketitle

 \addtocounter{footnote}{1} \footnotetext{Department of Mathematics, Saint Michael's College, One
  Winooski Park, Colchester, VT, 05458, USA and Department of
  Mathematics and Statistics, University of 
  Vermont, 16 Colchester Avenue, Burlington, VT,  05405, USA.
\texttt{jellis-monaghan@smcvt.edu}}
 \addtocounter{footnote}{1} \footnotetext{Instituto de Matem\'aticas, Universidad Nacional Aut\'onoma
  de M\'exico, Area de la Investigaci\'on Cient\'{\i}fica,
Circuito Exterior, C.U., Coyoac\'an, 04510 M\'exico D.F., M\'exico. 
\texttt{merino@matem.unam.mx}}

\section{Introduction}
\label{sec:introduction}

 We begin our exploration of graph polynomials and their applications
 with the Tutte polynomial, a renown tool for analyzing properties of
 graphs and networks.  This two-variable graph polynomial, due to
 W. T. Tutte~\cite{Tut47,Tut54,Tut67}, has the important universal
 property that essentially any multiplicative graph invariant with a
 deletion/contraction reduction must be an evaluation of it.  These
 deletion/contraction operations are natural reductions for many
 network models arising from a wide range of problems at the hearts of
 computer science, engineering, optimization, physics, and biology.

   In addition to surveying a selection of the Tutte polynomial's many
 properties and applications, we use the Tutte polynomial to showcase
 a variety of principles and techniques for graph polynomials in
 general. These include several ways in which a graph polynomial may
 be defined and methods for extracting combinatorial information and
 algebraic properties from a graph polynomial.  We also use the Tutte
 polynomial to demonstrate how graph polynomials may be both
 specialized and generalized, and how they can encode information
 relevant to physical applications.  


  We begin with the Tutte polynomial because it has a rich and
  well-developed theory, and thus it serves as an ideal model for
  exploring other graph polynomials in the next chapter, Graph
  Polynomials and Their Applications II: Interrelations and
  Interpretations. Furthermore, because of the Tutte polynomial's long
  history, extensive study, and its universality property, it is often
  a `point of contact' for research into other graph polynomials in
  that their study frequently includes exploring their relations to the
  Tutte polynomial.  These interrelationships will be a central theme
  of the following chapter.

In this chapter we give both recursive and generating function formulations of the Tutte polynomial, and state its universality in the form of a recipe theorem.  We give a number of properties and combinatorial interpretations for various evaluations of the Tutte polynomial.  We recover colorings, flows, orientations, network reliability, etc., and related polynomials as specializations of the Tutte polynomial.  We discuss the coefficients, zeros, and derivatives of the Tutte polynomial, and conclude with a brief discussion of computational complexity.


\section{Preliminary Notions}
\label{sec:preliminary}

 The graph terminology that we use is standard and generally follows
 Diestel~\cite{Die00}. Graphs may have loops and multiple
 edges.  For a graph $G$ we denote by
 $V(G)$ its set 
 of vertices and by $E(G)$ its set of edges. An
 oriented graph, $\vec{G}$, also called a digraph, has a direction 
 assigned to each edge. 

 \subsubsection{Basic Concepts}

 We first recall some of the notions of graph theory most used in this
 chapter.  Two graphs $G_1$ and $G_2$ are
 \emph{isomorphic}\index{isomorphic}, denoted
 $G_1\simeq G_2$, if there exists a bijection
 $\phi:V(G_1)\rightarrow 
 V(G_2)$ with $xy\in E(G_1)$ if and only if
 $\phi(x)\phi(y)\in 
 E(G_2)$. We denote by $\kappa(G)$ the number of
 connected components  
 of a graph $G$, and by $c(G)$ the number of \emph{non-trivial}\
 connected components, that is the number of connected components not
 counting isolated vertices. A graph is
 \emph{$k$-connected}\index{$k$-connected} if at least $k$ vertices
 must be removed to disconnect the
 graph.

 A \emph{cycle}\index{cycle} in a graph $G$ is a set of edges
 $e_1,\ldots, e_k$ such that, if $e_i=(v_i,w_i)$ for $1\leq i\leq k$,
 then $w_i=v_{i+1}$ for $1\leq i \leq k-1$; also $w_k=v_1$ and
 $v_i\neq v_j$ for $i\neq j$. A \emph{trail}\index{trail} is a path
 that may revisit a vertex, but
 not retrace an edge.  A \emph{circuit} \index{circuit} is
 a closed trail, and thus a cycle is just
 a circuit that does not revisit any
 vertices.  In the case of a digraph, the edges of a trail or circuit
 must be consistently oriented. 

 The dual notion of a cycle is that of \emph{cut}\index{cut} or
 \emph{cocycle}\index{cocycle}. If $\{V_1, V_2\}$ is a partition of
 the vertex set, and the set $C$, consisting of those edges with one end in $V_1$ and
 one end in $V_2$, is not empty, then $C$ is called a \emph{cut}. A
 cycle with one edge is called a \emph{loop}\index{loop} and a cocycle
 with one edge is called a
 \emph{cut-edge}\index{cut!-edge|see{bridge}} or
 \emph{bridge}\index{bridge}. We refer to an edge that is neither a
 loop nor a bridge as
 \emph{ordinary}\index{edge!ordinary}.
 
 A \emph{tree}\index{tree} is a connected graph without cycles. A
 \emph{forest}\index{forest} is a graph whose connected components are
 all trees. A subgraph $H$ of a graph $G$ is
 \emph{spanning}\index{graph!spanning} if
 $V(H)=V(G)$. Spanning trees 
 in connected graphs will play a fundamental role in the theory of the
 Tutte polynomial. Observe that a loop in a connected graph can be
 characterized as an edge that is in no spanning tree, while a bridge
 is an edge that is in every spanning
 tree.

 If $V' \subseteq V( G )$, then the
 \emph{induced 
 subgraph}\index{graph!induced} on 
 $V'$ has vertex set $V'$ and edge set those edges of $G$ with both
 endpoints in $V'$.  If $E' \subseteq E(G)$, then the
 \emph{spanning subgraph} induced by $E'$ has vertex set $V(G)$
 and edge set $E'$.

 \subsubsection{Deletion and Contraction}

 The operations of deletion and contraction of an edge are essential
 to the study of the Tutte polynomial. The graph obtained by deleting
 an edge \index{graph!deletion}$e\in E(G)$ is just
 $G\setminus e=(V,E\setminus e)$. The graph obtained by contracting an
 \index{graph!contraction} edge $e$ in $G$  results from identifying
 the endpoints of $e$ followed by removing $e$, and is denoted
 $G/e$. When $e$ is a loop, $G/e$ is the same as $G\setminus e$. It is
 not difficult to check that both deletion and contraction are
 commutative, and thus, for a subset of edges $A$, both $G\setminus A$
 and $G/A$ are well defined. Also, if $e\neq f$, then $G\setminus e/f$ and
 $G/f\setminus e$ are isomorphic; thus for disjoint subsets
 $A,A'\subseteq E(G)$, the graph $G\setminus A/A'$ is
 well-defined. A 
 graph $H$ isomorphic to $G\setminus A/A'$ for some choice of disjoint
 edge sets $A$ and $A'$ is called a \emph{minor}\index{graph!minor} of
 $G$. A class $\mathcal{G}$ of graphs  is \emph{minor closed}  if
 whenever $G$ is in  $\mathcal{G}$ then any 
 minor of $G$ is  also in the class.


 
  A \emph{graph invariant}\index{graph!invariant} is a function $f$ on
 the class of all graphs such that
\[
   f(G_1)=f(G_2)\ \text{whenever}\
   G_1\simeq G_2 .
\]
 A \emph{graph polynomial} \index{graph polynomial}  is a graph invariant where the image lies in
 some polynomial ring.

\subsubsection{The Rank and Nullity Functions for Graphs}

 To simplify notation, we typically identify a subset of edges $A$ of
 a graph $G$ with the spanning subgraph of $G$ that $A$ induces. Thus,
 for a fixed graph $G$ we have the following rank and nullity
 functions on the lattice of subsets of $E(G)$.

\begin{definition}\label{rank}\index{rank}\index{nullity}
 For $A\subseteq E(G)$, the rank and nullity of $A$,
 denoted 
 $r(A)$ and $n(A)$ respectively, are defined as
\[
    r(A)=|V(G)|-\kappa(A) \text{ and }
    n(A)=|A|-r(A) . 
\]
\end{definition}

Three special graphs are important. One is the rank 0 graph $L$
consisting of a single vertex with one loop edge, another is the
rank 1 graph $B$ consisting of two vertices with one bridge edge
between them, and the third one is the edgeless graph $E_1$ on $1$ vertex.

\subsubsection{Planar Graphs and Duality} 

 A graph is \emph{planar}\index{graph!planar} if it can be drawn in
 the plane without edges crossing, and it is a \emph{plane
 graph}\index{graph!plane} if it is so drawn in the plane. A drawing
 of a graph in the plane separates the plane into regions called
 faces. Every plane graph $G$ has a \emph{dual
 graph}\index{graph!dual}, $G^*$, formed by assigning a vertex of
 $G^*$ to each face of $G$ and joining two vertices of $G^*$ by $k$
 edges if and only if the corresponding faces of $G$ share $k$ edges
 in their boundaries. Notice that $G^*$ is always connected. If $G$ is
 connected, then $(G^*)^*=G$. If $G$ is planar, in principle it may
 have many plane duals, but when $G$ is 3-connected, all its plane
 duals are isomorphic. This is not the case when $G$ is only
 2-connected.

 There is a natural bijection between the edge set of a planar graph
 $G$ and the edge set of
 $G^*$, any one of its plane duals, so we can assume that $G$ and
 $G^*$ have the same edge set $E$. It is easy to check that
 $A\subseteq E$ is a spanning tree of $G$ if and only if $E\setminus
 A$ is a spanning tree of $G^*$. Thus, a planar graph and any of its
 plane duals have the same number of spanning trees. Furthermore, if
 $G$ is a planar graph with rank function $r$, and $G^*$ is any of
 its plane duals, then the rank function of $G^*$, denoted $r^*$,
 can be expressed as
\begin{equation}\label{eq:rankDuality}
   r^*(A)=|A|-r(E)+r(E\setminus A) .
\end{equation}

These observations reflect a deeper relation between $G$ and $G^*$
 that we will see captured by the Tutte polynomial at the end of
 Subsection~\ref{subsec:rank_generating}.

  



 \section{Defining the Tutte Polynomial}
\label{sec:III}

 Here we present several very different, but nevertheless equivalent,
 definitions of the Tutte polynomial.  The interplay among these
 different formulations is a source for many powerful tools developed
 to analyze the Tutte polynomial. Furthermore, each formulation lends
 itself to different proof techniques, for example induction with the
 linear recursion form and M\"obius inversion with the generating
 function form. These different formulations also are representative
 of some of the most common ways of defining any graph polynomial,
 although we will also see other methods in the next chapter.

While space prohibits including full proofs of the equivalence of
these various expressions for the Tutte polynomial, we note that there
are several approaches. One direct way is to specify the linear
recursion form as the definition of the Tutte polynomial and then use
induction on the number of edges to show that it is equivalent to
either the rank-nullity generating function or the spanning trees
expansion.  Showing that the linear recursion form is equivalent to
rank generating function form also establishes the essential fact that
it is well-defined, that is, independent of the order in which the
edges are deleted and contracted.  Another common approach is to
establish some definition of the Tutte polynomial, then prove from it
that the Tutte polynomial has the universality property discussed in
Section~\ref{sec:IV}. This universality property may then be applied
to show that some other function is equivalent to, or an evaluation
of, the Tutte polynomial.

 The spanning trees expansion formulation in
 Subsection~\ref{subsec:trees_expansion} was the approach originally
 used by Tutte to develop versions of this and similar polynomials.
 See Tutte~\cite{Tut47,Tut54,Tut67}.  A particularly lucid proof of
 the equivalence between rank-nullity generating function definition
 of Definition~\ref{def:rank_generating_expansion} and the spanning
 trees expansion definition of Definition~\ref{def:trees_expansion} can
 be found in~\cite{Bol98}. That the Tutte polynomial has a
 deletion-contraction reduction was shown by Tutte~\cite{Tut47, Tut48,
 Tut67} (also see Brylawski~\cite{Bry72}, and Oxley and
 Welsh~\cite{OW79}).

\subsection{Linear Recursion Definition}
\label{subsec:recursion}

 Broadly speaking, a linear recursion relation is a set of reduction
 rules together with an evaluation for the terminal forms.  The
 reduction rules rewrite a graph as a weighted (formal) sum of graphs
 that are in some way ``smaller'' or ``simpler'' than the original
 graph.  Furthermore, the reduction rules again apply to the newly
 generated simpler graphs, hence the recursion.  This recursion
 process eventually terminates in a well-defined set of ``most simple''
 graphs, which are no longer reducible by the reduction rules. These
 are then each identified with a monomial of independent variables to
 yield a polynomial.  It is essential to show that the reduction rules
 are independent of the order in which they are applied and that they
 do in fact terminate.  See Yetter~\cite{Yet90} for a more formal treatment.

 The Tutte polynomial may be defined by a linear recursion relation
 given by deleting and contracting ordinary edges. The ``most simple''
 terminal graphs are then just forests with loops.

\begin{definition}\label{def:recursive}
If $G=(V,E)$ is a graph, and $e$ is an ordinary edge, then
\begin{equation}\label{eq:deletion-contraction}
   T(G;x,y) = T(G\setminus e; x, y) + T(G/e; x, y).
\end{equation} 
 Otherwise, $G$ consists of $i$ bridges and $j$ loops and
\begin{equation} \label{end_recursion}
 T(G; x, y) = x^{i}y^{j}.
\end{equation}
\end{definition}

 In other words, $T$ may be calculated recursively by specifying an
 ordering of the edges and repeatedly
 applying~(\ref{eq:deletion-contraction}).  Remarkably, the Tutte
 polynomial is well defined in that the polynomial resulting from this
 recursive process is independent of the order in which the edges are
 chosen. One way to prove this is by showing that this definition is
 equivalent to the rank generating form we will see in
 Definition~\ref{def:rank_generating_expansion}. A proof can be found
 in~\cite{BO92} for example.

\begin{figure}[hbtp]
 \begin{center}
 \includegraphics[scale=0.25]{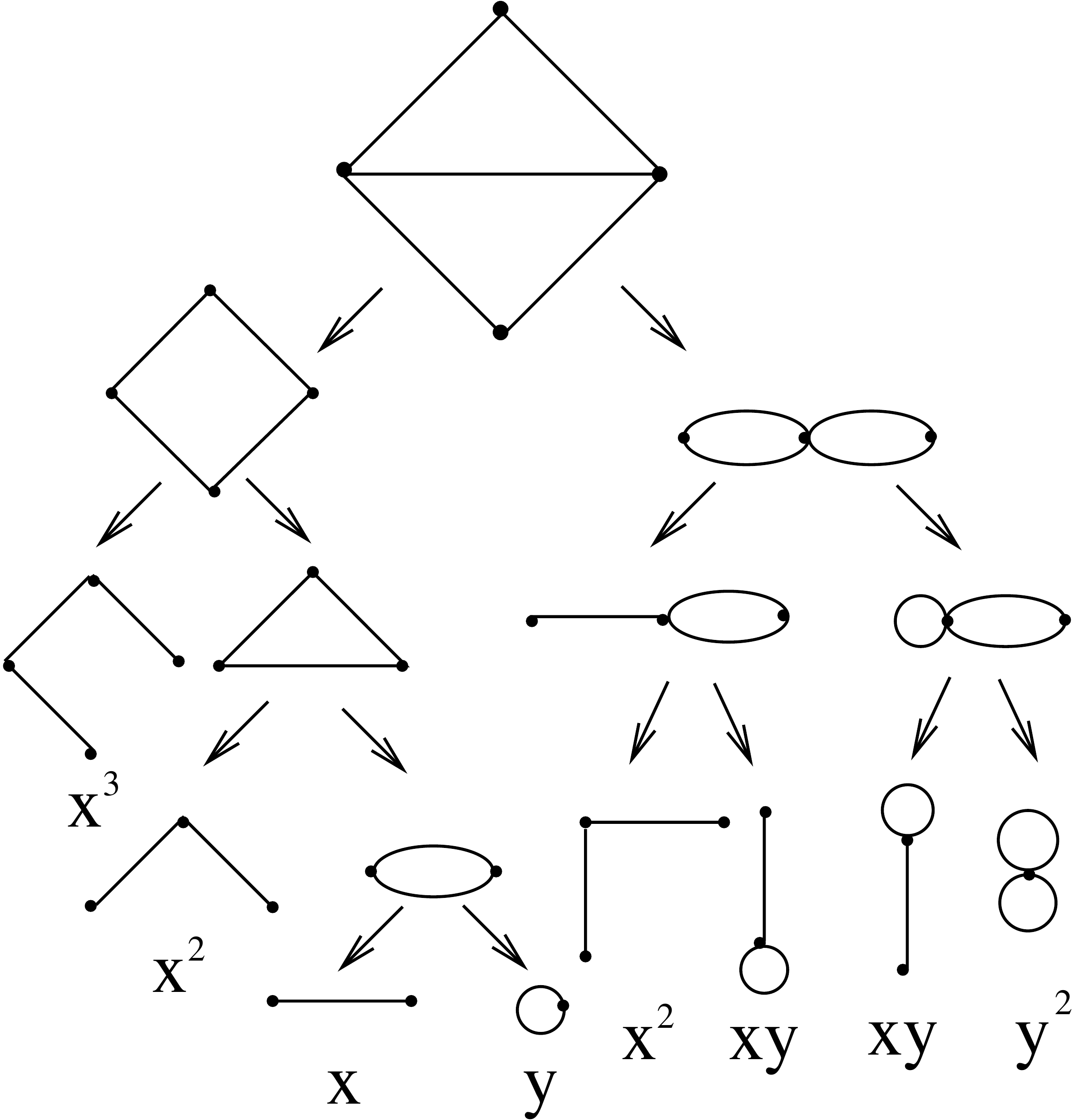}
 \caption{An example of computing the Tutte polynomial
   recursively}\label{Fig:recursive}
 \end{center}
 \end{figure}

 Figure~\ref{Fig:recursive} gives a small example of computing $T$ using
 (\ref{eq:deletion-contraction}) and (\ref{end_recursion}) for $K_4$
 minus one edge. By adding the monomials at the bottom of
 Figure~\ref{Fig:recursive} we find that $T(G; x, y)= x^3+2 x^2+ x
 + 2 x y + y+ y^2$.

Recall that a \emph{one-point join} \index{graph!one-point join} $G*H$
of two graphs $G$ and $H$ is formed by identifying a vertex $u$ of $G$
and a vertex $v$ of $H$ into a single vertex $w$ (necessarily a cut
vertex) of $G*H$. Also, $G \cup H$ is the disjoint union  of $G$ and
$H$. 

  \begin{proposition}\label{product-formula}
  If $G$ and $H$ are graphs then
 \[
    T(G \cup H)=T(G)\,T(H)
    \textit{ and }
    T(G*H)=T(G)\,T(H).
 \]
  \end{proposition}

This follows readily from Definition~\ref{def:recursive} by induction
on the number of ordinary edges in $G*H$ or $G \cup H$.
 
\subsection{Rank-Nullity Generating Function Definition}
\label{subsec:rank_generating}
 A generating function can often be thought of as a (possible
 infinite) polynomial whose coefficients count structures that are
 encoded by the exponents of the variables.  Because generating
 functions count, which is at the very heart of enumerative
 combinatorics, there is extensive literature on them. The two
 volumes~\cite{Sta96} and~\cite{Sta99} are an excellent resource.  In
 the case of the Tutte polynomial, there are several different
 generating function formulations, each of which has its advantages.
 We give one here and another in
 Subsection~\ref{subsec:trees_expansion}, with a variation in
 Subsection~\ref{subsec:coefficient_relations}, and refer the reader
 to~\cite{BO92} for additional forms.

\begin{definition}\label{def:rank_generating_expansion}
 If $G=(V,E)$ is a graph, then the Tutte polynomial of $G$, $T(G;$
 $x$, $y)$, has the following expansion
 \begin{equation}\label{eq:rank_expansion}
     T(G; x, y) = \sum_{A\subseteq E} (x-1)^{r(E)-r(A)}
                                      (y-1)^{n(A)}.
 \end{equation}
 \end{definition}

 The advantage of a generating function formulation is that it
 facilitates counting.  For example, interpretations for several
 evaluations of the Tutte polynomial given in Section~\ref{sec:V}
 follow immediately from
 Definition~\ref{def:rank_generating_expansion}.

We can also deduce the following pleasing property of the Tutte
polynomial.

\begin{proposition}\label{planar-dual}
  If $G$ is a planar graph with dual $G^*$ then
   \begin{equation}\label{eq:6::duality}
        T(G;x,y) = T(G^*;y,x).
   \end{equation}
\end{proposition}

 This result follows from routine checking using
 Definition~\ref{def:rank_generating_expansion} and
 (\ref{eq:rankDuality}).

\subsection{Spanning Trees Expansion Definition}
\label{subsec:trees_expansion}

 We need to develop a little terminology before presenting the 
 spanning trees definition of the Tutte polynomial. First, given a
 spanning tree $S$ and an edge $e\not\in S$, there is a cycle defined
 by $e$, namely the unique cycle in $S\cup e$. Similarly, for an edge
 $f\in S$, there is a cut defined by $f$, namely the set of edges
 $C$ such that if $´f'\in C$, then $(S-f)\cup f'$ is a spanning tree.

 Assume there is a fixed ordering $\prec$ on the edges of $G$, say
 $E=\{e_1,\ldots,e_m\}$, where $e_i\prec e_j$ if $i<j$. Given a fixed
 tree $S$, an edge $f$ is called {\it internally active}
 \index{active edge!internally} if $f\in S$ and it is the smallest
 edge in the
 cut defined by $f$. Dually, an edge $e$ is {\it externally
 active}\index{active edge!externally} if $e\not\in S$ and it is the
 smallest edge in the cycle defined by $e$. The internal activity of
 $S$ is the number of its internally active edges and its external
 activity is the number of externally active edges.  With this, we
 have the following definition of  the Tutte polynomial.

\begin{definition}\label{def:trees_expansion}
 If $G$ is a graph with a total order on its edge set, then
\begin{equation}
   T(G;x,y)=\sum_{i,j}t_{ij}x^iy^j,
\end{equation}
 where $t_{ij}$ is the number of spanning trees with internal activity
 $i$ and  external activity $j$.
\end{definition} 

Two important observations follow immediately from the equivalence of
Definitions~\ref{def:rank_generating_expansion}
and~\ref{def:trees_expansion}.  One is that the the terms $t_{ij}$ in 
Definition~\ref{def:trees_expansion} are independent of the total
order used in the edge set, since there is no ordering of the edges in
Definition~\ref{def:rank_generating_expansion}. The other is that the
coefficients in Definition~\ref{def:rank_generating_expansion} must be
non-negative since the coefficients in
Definition~\ref{def:trees_expansion} clearly are.

 \section{Universality of the Tutte Polynomial}
\label{sec:IV}

The universality property discussed here is one of the most powerful
aspects of the Tutte polynomial.  It says that essentially any graph 
invariant that is multiplicative on disjoint unions and one-point
joins of graphs and that has a deletion/contraction reduction must be
an evaluation of the Tutte polynomial. We will see several
applications of this theorem throughout the rest of this chapter and
in the next chapter as well. Various generalizations of the Tutte
polynomial are careful to retain this essential property, and
analogous universality properties are sought in the context of other
graph polynomials.

 \begin{definition}\label{T-G}\index{Tutte-Gr\"othendieck invariant}
 Let $\mathcal{G}$ be a minor closed class of graphs.
  A  graph invariant $f$ from $\mathcal{G}$ to a
 commutative ring $\mathcal{R}$ with unity  is called a
 \emph{generalized Tutte-Gr\"othendieck invariant},
  or \emph{T-G invariant}, if $f(E_1)$ is the unity of $R$, if
 there exist fixed  elements $a,b \in R$ such that for every
 graph $G\in \mathcal{G}$ and every ordinary edge $e\in G$,  then
 \begin{equation} \label{eq:TG_invariant_del_con}
     f(G)= a f(G\setminus e)+
     bf(G/e) ;
 \end{equation} 
  and if for every $G,H \in \mathcal{G}$, whenever $G \cup H$ or $G*H$
  is in $\mathcal{G}$, then
 \begin{equation} \label{eq:TG_invariant_prod}
    f(G \cup H)=f(G)f(H)
    \textit{ and }
    f(G*H)=f(G)f(H).
 \end{equation}
\end{definition}

Thus, the Tutte polynomial is a T-G invariant, and in fact, since the
following two results give both universal and unique extension
properties, it is essentially the \emph{only} T-G invariant, in that
any other must be an evaluation of it. Theorem~\ref{universal} is
known as a recipe theorem since it specifies  
how to recover a T-G invariant as an evaluation of the Tutte polynomial.

\begin{theorem}\label{universal}
\index{Tutte polynomial!universal property}
Let $\mathcal{G}$ be a minor closed class of graphs, let $R$ be a
commutative ring with unity, and let $f:\mathcal{G} \to
R$.  If there exists $a,b \in R$ such that $f$ is a T-G
invariant, then
\begin{equation}\label{recipe}
f(G) =   
          a^{|E(G)|- r(E(G))}
          b^{r(E(G))}
          T\left({G;\frac{{x_0 }}{b},
              \frac{{y_0 }}{a}} \right) ,
\end{equation}
 where $f(B) = x_0 $ and 
$f\left( L \right) =y_0 $.
\end{theorem}

Furthermore, we have the following unique extension property, which
says that if we specify any four elements $a, b, x_0, y_0 \in R$, then
there is a unique well-defined T-G invariant on these four elements.

\begin{theorem}\label{unique}
\index{Tutte polynomial!uniqueness property} Let $\mathcal{G}$ be a
  minor closed class of graphs, let $R$ be a commutative ring with
  unity, and let $a, b, x_0, y_0 \in R$. 
  Then  there is a unique T-G  invariant $f:\mathcal{G} \to R$
  satisfying Definition ~\ref{T-G}  with $f(B)=x_0$ and
  $f(L)=y_0$. Furthermore, this function
  $f$ is given by
 \begin{equation}\label{recipe_2}
      f(G)= a^{|E(G)|- r(E(G))} 
          b^{r(E(G))}
               T\left( G;\frac{x_0}{b}, \frac{y_0}{a}\right) .
  \end{equation}
\end{theorem}

If $a$ or $b$ are not units of $R$, then (\ref{recipe}) and (\ref{recipe_2})
are interpreted 
to mean using  expansion (\ref{eq:rank_expansion}) of
Definition~\ref{def:rank_generating_expansion}, and cancelling before
evaluating.  

 These results can be proved by induction on the number of ordinary
 edges from the deletion/contraction definition of the Tutte
 polynomial.  See, for example, Brylawski~\cite{Bry72}, Oxley and
 Welsh~\cite{OW79}, Brylawski and Oxley~\cite{BO92}, 
 Welsh~\cite{Wel93}, and Bollob\'as~\cite{Bol98} for detailed
 discussions of these theorems and  their consequences.

Examples applying this important theorem may found throughout
Section~\ref{specializations}, where it may be used to show that all
of the graph polynomials surveyed there are evaluations of the Tutte
polynomial.


 \section{Combinatorial Interpretations of Some Evaluations}
\label{sec:V}

A graph polynomial encodes information about a graph.  The challenge
is in extracting combinatorially useful information from this
algebraic object.  A number of successful techniques have evolved for
meeting this challenge, and we use the Tutte polynomial to showcase
some of them while simultaneously demonstrating the richness of the
information encoded by the Tutte polynomial.

\subsection{Spanning Subgraphs}
\label{subsec:trees}

 Spanning subgraphs, and in particular spanning trees, play a
 fundamental role in the theory of Tutte polynomials as we have
 already seen in Definition~\ref{def:trees_expansion}. This is also
 reflected in the most readily attainable interpretations for
 evaluations of the Tutte polynomial, which enumerate various spanning
 subgraphs.  We begin with these here, writing $\tau(G)$ for the
 number of spanning trees of a connected graph $G$.

\begin{theorem}\label{trivial_interpretations}
 If $G=(V,E)$ is a connected graph then:

 \begin{enumerate}
 \item $T(G;1,1)$ equals $\tau(G)$.

 \item $T(G;2,1)$ equals the number of spanning forests of $G$.
 \label{span2-1}

 \item $T(G;1,2)$ equals the number of spanning connected
       subgraphs of $G$.
     
\item $T(G;2,2)$ equals $2^{|E|}$.
 \end{enumerate}
 
 \end{theorem}
\begin{proof}
 To illustrate common proof techniques, we give two short proofs of
 the first statement. The remaining statements
 may be proved similarly. When $x=y=1$, the non-vanishing terms in
 the rank-nullity expansion (\ref{eq:rank_expansion}) are
 $A\subseteq E$ such that $r(E)=r(A)$ and $|A|=r(A)$.
 That $r(E)=r(A)$ implies that $A$ has the same number of
 connected components as $G$, namely one, so $(V,A)$ is connected.
 Then $|A|=r(A)$ implies that $|A|=|V|-1$, so $A$ must be a tree,
 and hence a spanning tree.

 Alternatively, we can use  Theorem~\ref{universal}. Let $\tau'(G)$
 be the number of maximal spanning forests in a general (not
 necessarily connected) graph $G$. We
 prove that $T(G;1,1)=\tau'(G)$. If $G$ is connected, we have that
 $T(G;1,1)=\tau'(G)=\tau(G)$. First note that the number of maximal
 spanning  forests has a 
 deletion-contraction reduction for ordinary edges, that is, if $G$ is
 a graph and $e$ is an ordinary edge of $G$, then
 $\tau'(G)=\tau'(G\setminus e)+\tau'(G/e)$.  This follows because the
 maximal spanning forests of $G$ can be partitioned into the maximal
 spanning forests that 
 do not contain $e$ and those that do contain $e$. The former are  the
 maximal spanning forests of $G\setminus e$ 
 and the latter are in one-to-one correspondence with the maximal spanning
 forests of $G/e$. 

The result then follows immediately from Theorem~\ref{universal} with
$a=b=x_0=y_0=1$. \qed
\end{proof}

Computing the number of spanning trees of a graph is easy in that
there are polynomial time algorithms to do it. One of these involves
a determinant. Recall that the \emph{Laplacian matrix} $L$ of an
graph $G$ with vertices $v_1, \dots, v_n$ is the $n\times n$-matrix
defined by 
\begin{equation}\label{Laplacian}
\index{Laplacian matrix} L_{ij}=\begin{cases} \operatorname{deg}(i)&
  \text{if 
$i=j$}\\ -r& \text{if $r$ is the number of edges between vertices $i$
and $j$}.
        \end{cases}
\end{equation}

\begin{theorem}\label{cofactors} If $G$ is a connected graph with
  Laplacian $L$, then 
\begin{equation}
   T(G;1,1)=\tau(G)=\operatorname{Det}(L'),
\end{equation}
 where $L'$ is any cofactor of $L$.
\end{theorem}
A proof of this can be found in~\cite{AZ01} using the Binet-Cauchy
 formula and that $L=DD^t$, where $D$ is the incidence matrix of (an
 orientation of) $G$.

This result not only provides an interpretation of the Tutte
polynomial at $(1,1)$ in terms of the incidence matrix of a graph, but
also proves that $(1,1)$ is one of the (very few, as we will see in
Section~\ref{sec:complexity}) points where the Tutte polynomial can be
computed in polynomial time.

\subsection{The Tutte Polynomial at $y=x$}
\label{subsec:x=y}

Combinatorial interpretations are known for the Tutte polynomial at
all integer values along the line $y=x$.  In addition to those for
$T(G,1,1)$ and $T(G,2,2)$ previously given, we have the following
interpretation for $T(G;-1,-1)$ due to Read and
Rosenstiehl~\cite{RR78}.  We will also see alternative interpretations
for $T(G;-1, -1)$ and $T(G;3,3)$ in Subsection~\ref{orientations}.

  The incidence matrix $D$ of a graph $G$ defines a vector
  space over 
  $\mathbb{Z}_2$, called the cycle space \index{graph!cycle space}
  $\mathcal{C}$. The bicycle
  space\index{graph!bicycle space} $\mathcal{B}$ is then just
  $\mathcal{C}\cap \mathcal{C}^{\perp}$.

\begin{theorem} \label{bicycle}
  $T(G;-1,-1)=(-1)^{\left|E\right|}(-2)^{\dim(\mathcal{B})}$. 
\end{theorem}

 
One method of extracting information from a graph polynomial is via
 its relation to some other graph invariant.  The following
 interpretations for the Tutte polynomial of a planar graph along the
 line $x=y$ derive from its relation to the Martin
 polynomial\index{Martin polynomial}\index{polynomial!Martin}
 \cite{Mar77}, a one variable graph polynomial that we will discuss
 in further in the next chapter.

\begin{figure}[hbtp]
 \begin{center}
 \includegraphics[scale=0.5]{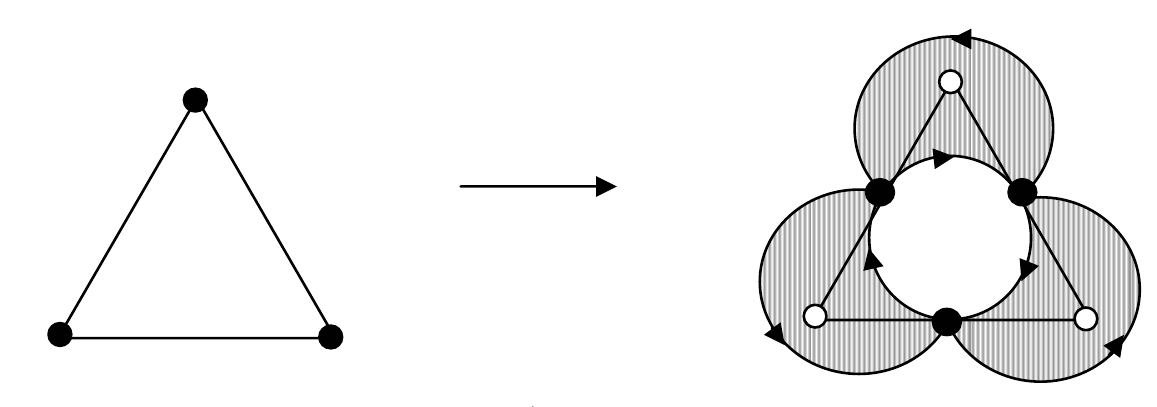}
 \caption{On the left hand side we have a planar graph $G$. On the right
   hand side we have $\vec{G}_m$ with the vertex faces colored black,
   oriented 
   so that black faces are to the left of each edge}\label{Fig:medial} 
 \end{center}
 \end{figure}

 We first recall that the \emph{medial graph}\index{graph!medial} of a
 connected  planar graph $G$ 
 is constructed by placing a vertex on each edge of $G$ and drawing
 edges around the faces of $G$.  The faces of this medial graph are
 colored black or white, depending on whether they contain or do not
 contain, respectively, a vertex of the original graph $G$.  This face
 two-colors the medial graph.  The edges of the medial graph are then
 directed so that the black face is on the left.  We refer to this as
 the directed medial graph of $G$ and denote it by $\vec{G}_m$. An
 example is given in Fig.~\ref{Fig:medial}.
 
 Martin~\cite{Mar77} showed that, for a planar graph $G$, the relation
 between the Martin polynomial and Tutte polynomial is
 ${m}(\vec{G}_m;x)=T(G;x,x)$.  Evaluations for the
 Martin 
 polynomial in~\cite{E-M04} then give the following interpretations of
 the Tutte polynomial.

 Let ${D}_{n}(\vec{G}_m)=\{(D_1,\ldots,D_n)\}$, where
 $(D_1,\ldots,D_n)$ 
 is an ordered partition of $E(\vec{G}_m)$ into $n$
 subsets such that 
 $G$ restricted to $D_i$ is 2-regular and consistently oriented for all $i$.

 \begin{theorem}
   Let $G$ be a planar graph with oriented medial graph $\vec{G}_m$.
  Then, for $n$ a positive integer,
   \[
    (-n)^{c(G)}T(G;1-n,1-n)=\sum_{{D}_n(\vec{G}_m)}
    (-1)^{\sum_{i=1}^{n}c(D_i)}.
     \]
 \end{theorem}

 \begin{theorem}
   Let $G$ be a planar graph with oriented medial graph $\vec{G}_m$.
  Then, for $n$ a positive integer,
   \[
    n^{c(G)}T(G;1+n,1+n)=\sum 2^{\mu(\phi)},
     \]
  where the sum is over all edge colorings $\phi$ of $\vec{G}_m$ with
  $n$ colors so that each (possibly empty) set of monochromatic edges
  forms an Eulerian digraph, and where $\mu(\phi)$ is the number of
  monochromatic vertices in the coloring $\phi$.
 \end{theorem}

\subsection{Orientations and Score Vectors}\label{orientations}

  The combinatorial interpretations of the Tutte polynomial in
  Theorem~\ref{trivial_interpretations} are given in terms of the
  number of certain subgraphs of the graph $G$. However, they can also be
  given in terms of orientations of the graph and its score
  vectors. Given a graph $G=(V,E)$, an orientation of
  $G$\index{graph!orientation} may be obtained by directing every edge
  from either of its ends to the other.  From this follows that
  $T(G;2,2)$ equals the number of possible orientations of $G$.

  The \emph{score vector} \index{score vector} of an orientation
  $\vec{G}$ is the vector ($s_1$, $s_2$, \ldots,$s_n$) such that
  vertex $i$ has outdegree $s_i$ in the orientation.  In the following
  theorem we gather several similar results about the Tutte polynomial
  and orientations of a graph.

\begin{theorem}\label{orientation_results}
 Let $G=(V,E)$ be a connected graph with Tutte polynomial $T(G;$
 $x$, $y)$. Then

 \begin{enumerate}
 \item $T(G;2,0)$ equals the number of acyclic orientations
 \label{2-0} \index{orientation!acyclic} of $G$, that is orientations
 without oriented cycles.

 \item $T(G;0,2)$ equals the number of totally cyclic
 orientations \index{orientation!totally cyclic} of $G$, that is
 orientations in which every arc is in a directed cycle\label{0-2}.

 \item $T(G;1,0)$ equals the number of acyclic orientations with
 exactly one predefined source $v$\label{1-0}.

 \item $T(G;2,1)$ equals the number of score vectors of
       orientations of $G$\label{2-1}.
     

 \end{enumerate}
\end{theorem}

Item~\ref{2-0} was first proved by Stanley in~\cite{Sta73} by using
 the Ehrhart polynomial; a proof using the universality of the Tutte
 polynomial is given by Brylawski and Oxley in~\cite{BO92}. For
 Items~\ref{0-2} and~\ref{1-0}, see Green and Zaslavsky~\cite{GZ83}
 and Las Vergnas~\cite{Las77}. The former also gives the number of
 strongly connected orientations of $G$, and note that the latter is
 independent of the choice of source vertex $v$. Item~\ref{2-1} was
 first proved by Stanley in \cite{Sta80}, with a bijective proof given
 by Kleitman and Winston in \cite{KW81}, and a proof using
 Theorem~\ref{universal} by Brylawski and Oxley in
 \cite{BO92}. Comparing Item~\ref{2-1} with
 Theorem~\ref{trivial_interpretations}, Item~\ref{span2-1}, shows that
 the number of score vectors equals the number of spanning forests of
 $G$.

Some other evaluations of the Tutte polynomial can also be interpreted
 in terms of orientations. Recall that an anticircuit
 \index{anticircuit} in a digraph is a closed trail so that the
 directions of the edges alternate as the trail passes through any
 vertex of degree greater than 2.  Note that in a 4-regular Eulerian
 digraph such as  $\vec{G}_m$, the set of anticircuits can be found by
 pairing the two incoming edges and the two outgoing edges at each
 vertex.

 Two surprising results from Las Vergnas~\cite{Las88} and
 Martin~\cite{Mar78} are the following:

\begin{theorem}  Let $G$ be a connected planar graph.  Then

   \begin{equation}\label{minusone}
 T(G;-1,-1)=(-1)^{|E(G)|}
                (-2)^{{a}(\vec{G}_m)-1}  
   \end{equation}

and
    \begin{equation}\label{three}
    T(G;3,3)=K 2^{{a}(\vec{G}_m)-1},
    \end{equation}
 where ${a}(\vec{G}_m)$ is the number of anticircuits in
 the directed medial graph of G and $K$ is some odd integer.
\end{theorem}

Comparing (\ref{minusone}) to Theorem \ref{bicycle} gives the following
corollary.

\begin{corollary}
If $G$ be a connected planar graph, then the dimension of the bicycle
 space is ${a}(\vec{G}_m)-1$.
\end{corollary}

 \section{Some Specializations}
\label{specializations}
Here we illustrate the wide range of applicability of the Tutte
 polynomial while demonstrating some proof techniques for showing that
 a graph invariant is related to the Tutte polynomial.  The advantage
 of recognizing an application-driven function as a specialization of
 the Tutte polynomial is that the large body of 
 knowledge about the Tutte polynomial is then available to inform the
 desired application. We say a graph polynomial is a specialization of
 the Tutte polynomial if it may be recovered from the Tutte polynomial
 by some substitution for $x$ and $y$, with possibly some prefactor.  

For various substitutions along different algebraic curves in $x$ and
 $y$, the Tutte polynomial has interpretations as the generating
 function of combinatorial quantities or numerical invariants
 associated with a graph. Some of these were considered long before
 the development of the Tutte polynomial, and others were discovered
 to be unexpectedly related to the Tutte polynomial.  We survey six of
 the more well known of these application-driven generating functions.

\subsection{The Chromatic Polynomial}

The chromatic polynomial, because of its theoretical and applied
 importance, has generated a large body of work. Chia \cite{Chi97}
 provides an extensive bibliography on the chromatic polynomial, and
 Dong, Koh, and Teo \cite{DKT05} give a comprehensive treatment.
 
For positive integer $\lambda$, a $\lambda$-coloring
\index{coloring@$\lambda$-coloring} of a graph $G$
is a mapping of $V(G)$ into the set
$[\lambda]=\{1,2,\ldots,\lambda\}$. Thus there are exactly
$\lambda^{n}$ colorings for a graph on $n$ vertices. If $\phi$ is a
$\lambda$-coloring such that $\phi(i)\neq \phi(j)$ for all $ij\in E$,
then $\phi$ is called a \emph{proper} (or \emph{admissible})
coloring.

 We wish to find the number, $\chi(G;\lambda)$, of admissible
$\lambda$-colorings of a graph $G$. As noted by Whitney~\cite{Whi32},
the 4-color theorem can be formulated in this general setting as
follows: If $G$ is planar graph, then $\chi(G;4)>0$.

 The following theorem is due to G. D. Birkhoff in~\cite{Bir12} and
 independently by H. Whitney in~\cite{Whi32}. We sketch the second
 proof.
\begin{theorem}\label{TheoWhitney} 
 If $G=(V,E)$ is a graph, then
 \begin{equation}\label{whitney}
    \chi(G;\lambda)= \sum_{A\subseteq E}{(-1)^{|A|}
    \lambda^{\kappa(A)}}.
 \end{equation}
 \end{theorem}
\begin{proof}
 Let $P_{ij}$ be the set of $\lambda$-colorings such that vertices
 $i$ and $j$ receive the same color. Let $\bar{P}_{ij}$ be the
 complement of $P_{ij}$ in the set of $\lambda$-colorings. Then, the
 value $\chi(G;\lambda)$ can be computed using the inclusion-exclusion
 principle.
 \begin{eqnarray} \label{inc_exc}
   \chi(G;\lambda)&=& \left|\bigcap_{ij\in E}
            {\bar{P}_{ij}}\right|\nonumber\\ 
          &=& \lambda^{n} - \sum\limits_{ij \in E}
              {\left| {P_{ij} } \right|} +
              \sum\limits_{\begin{subarray}{l} ij,kl \in E \\ ij
              \ne kl \end{subarray}} 
              {\left| {P_{ij}  \cap P_{kl} } \right|}-\cr 
          & \cdots & 
              +(-1)^{|E|}\left|\bigcap_{ij\in E} {P_{ij}}\right|. 
 \end{eqnarray}

 Every term is of the form 
 $\left|\cap_{ij \in A} {P_{ij}}\right|$ 
 for some $A \subseteq E$, and hence corresponds to the
  subgraph $(V,A)$, where $A$ is the set of edges given by the indices
  of the $P_{ij}$'s. Thus, the cardinality of this set is the number
  of $\lambda$-colorings that have a constant value on each of the
  connected components of $(V,A)$, that is, $\lambda^{\kappa(A)}$. The
  sum on the right-hand side of (\ref{whitney}) is then precisely
  (\ref{inc_exc}).  \qed\end{proof}

 Thus, $\chi(G;\lambda)$ is a polynomial on $\lambda$ and it is called
 the \emph{chromatic  polynomial} of $G$.
 Some easily seen properties of $\chi(G;\lambda)$, which can be found
 in Read's seminal work~\cite{Rea68} on the chromatic polynomial, are
 the following:

\begin{proposition}\label{chrom_props} If $G$ is a graph with
 chromatic polynomial $\chi(G;\lambda)$, then: 
 \begin{enumerate}
 \item If $G$ has no edges, then
 $\chi(G;\lambda)=\lambda^{n}$\label{chrom terminal}. 
 \item If $G$ has a loop, then $\chi(G;\lambda)=0$, for all $\lambda$.
 \item $\chi(K_{n};\lambda)=\lambda(\lambda-1)\cdots(\lambda-n+1)$.
 \item \label{del_con_col} If $e$ is any edge of $G$, then
  \[ 
   \chi(G;\lambda)= \chi(G\setminus e;\lambda)-
   \chi(G/e;\lambda)\label{chrom recursion}. 
 \]
 \end{enumerate}
\end{proposition}

 Note that Items~\ref{chrom terminal} and~\ref{del_con_col} give a
 recursive alternative definition of the chromatic polynomial.

 Also in Read~\cite{Rea68} is the following not so trivial, but not
 difficult to prove, property of the chromatic polynomial. 

\begin{theorem}\label{ChrJoint}
  If $G$ is the union of two vertex set induced subgraphs $H_1$ and
  $H_2$ such that the intersection $H_1\cap H_2$ is a vertex set
  induced subgraph isomorphic to $K_p$, then
\[
  \chi(G;\lambda)=
   \frac{\chi(H_1;\lambda)\chi(H_2;\lambda)}{\chi(K_{p};\lambda)}.
\]
\end{theorem}

 Thus, although Proposition~\ref{chrom_props} Item~\ref{del_con_col}
 suggests that the chromatic polynomial might be a T-G invariant, by
 Theorem~\ref{ChrJoint}, it  is \emph{not}
 multiplicative on the one point join of two graphs.  However, as is
 frequently the case, this can be addressed by a simple multiplier; it
 is easy to check that $\lambda^{-\kappa(G)}\chi(G; \lambda)$ is a T-G
 invariant. The relation between the Tutte and chromatic polynomials
 may then be found by  applying
 Theorem~\ref{universal}  with the help of
 Proposition~\ref{chrom_props}, Item~\ref{del_con_col}. We give an
 alternative proof of this relationship deriving from
 Theorem~\ref{TheoWhitney}.

\begin{theorem}
 If $G=(V,E)$ is a graph, then
 \[
   \chi(G; \lambda) = (-1)^{r(E)}\lambda^{\kappa(G)}T(G;
   1-\lambda, 0).
 \]
\end{theorem}
\begin{proof}
 Since $r(E)-r(A)=\kappa(A)-\kappa(G)$ we have that
  \begin{eqnarray*}
    \chi(G;\lambda) &=& \sum_{A\subseteq E}{(-1)^{|A|}
                        \lambda^{\kappa(A)}}\\
                        &=&(-1)^{r(E)}\lambda^{\kappa(G)}
                        \sum_{A\subseteq
                        E}{(-1)^{|A|-r(A)}(-\lambda)^{r(E)-r(A)}}\\ 
                        &=&(-1)^{r(E)}\lambda^{\kappa(G)}T(G;
                        1-\lambda, 0),
  \end{eqnarray*}
with the last equality following from
Definition~\ref{def:rank_generating_expansion}.  \qed\end{proof}

\subsection{The Bad Coloring Polynomial}

 One way to generalize the chromatic polynomial is to count \emph{all}
 possible colorings of the graph $G$, not just proper colorings. In
 order to differentiate between proper and improper colorings, we
 keep track of the edges between vertices of the same color, calling
 them \emph{bad edges}. This leads to the \emph{bad coloring polynomial}.

\begin{definition}\label{bad coloring}\index{bad coloring polynomial}
 \index{polynomial!bad coloring}

The \emph{bad coloring polynomial} is 
the generating function
 \[
    B(G; \lambda, t)= \sum_{j}{b_j(G; \lambda) t^j},
 \]
 where $b_{j}(G; \lambda)$ is the number of $\lambda$-colorings of
 $G$ with exactly $j$ bad edges. 
\end{definition}


Now consider $B(G; \lambda, t+1)$, which can be written as
 \begin{equation}\label{eq:gen_colouring}
       B(G; \lambda, t+1)=\sum_{\phi:V\rightarrow
         [\lambda]}{(1+t)^{|{b}(\phi)|}},
 \end{equation}
 where ${b}(\phi)$ is the set of bad edges in the $\lambda$-coloring
 $\phi$. With this last expression is again easy to get the relation
 to the Tutte polynomial using the following derivation of
 S. D. Noble (private communication).
\begin{theorem} 
 For a graph $G=(V,E)$ we have that
 \[
   B(G; \lambda, t+1) = t^{r(E)}\lambda^{\kappa(G)}
   T\left(G;\frac{\lambda+t}{t} ,1+t\right).
 \]
\end{theorem}
\begin{proof}
  \begin{eqnarray*}
     B(G; \lambda, t+1)&=&\sum_{\phi:V\rightarrow
                          [\lambda]}{(1+t)^{|{b}(\phi)|}}\\
                          &=&\sum_{\phi:V\rightarrow [\lambda]}
                          \sum_{A\subseteq {b}(\phi)}
                          t^{|A|}\\ 
                          &=&\sum_{A\subseteq
                          E}\sum_{\substack{\phi:V\rightarrow
                          [\lambda]\\ A\subseteq
                          {b}(\phi)}} t^{|A|}\\ 
                          &=& \sum_{A\subseteq E} t^{|A|}
                          \lambda^{\kappa(A)}.
  \end{eqnarray*}

Thus,
  \begin{eqnarray*}
    B(G; \lambda, t+1) 
           &=& \sum_{A\subseteq E}{t^{|A|}\lambda^{\kappa(A)}}\\
           &=&t^{r(E)}\lambda^{\kappa(G)}
              \sum_{A\subseteq E}{t^{|A|-r(A)}
              \left(\frac{\lambda}{t}\right)^{r(E)-r(A)}}\\
           &=&t^{r(E)}\lambda^{\kappa(G)}
              T\left(G;1+\frac{\lambda}{t},t+1\right). 
    \end{eqnarray*}
  \qed
\end{proof}

 Again, the above result could also be obtained from the universal
 property of the Tutte polynomial given in Theorem~\ref{universal} by
 applying it to $ \bar{B}(G; \lambda, t)=\lambda^{-\kappa(G)}
 B(G; \lambda, t)$ and verifying that

 \begin{enumerate}
 \item $\bar{B}(G; \lambda, t) = \bar{B}(G\setminus e; \lambda,
                     t) + (t-1) \bar{B}(G/e; \lambda, t)$, if $e$
                     is an ordinary edge.

 \item\label{6:15} $\bar{B}(G; \lambda, t) = t\bar{B}(G\setminus
 e; \lambda, t)$, if $e$ is a loop.

 \item\label{6:16} $\bar{B}(G; \lambda, t) = (t+\lambda
 -1)\bar{B}(G/e; \lambda, t)$, if $e$ is a bridge.
 \end{enumerate}

\subsection{The Flow Polynomial}
\label{flow_polynomial}

 The dual notion to a proper $\lambda$-coloring is a nowhere zero
 $\lambda$-flow. A standard resource for the
 material in this subsection is Zhang~\cite{Zha97}, while
 Jaeger~\cite{Jae88} gives  a good  survey.

Let $G$ be a graph with an arbitrary but fixed orientation, and let
 $H$ be an Abelian group with $0$ as its identity element. An
 $H$-flow\index{flow!H-flow@$H$-flow} is a mapping $\phi$ of the
 oriented edges $\vec{E}(G)$ into the elements of the
 group $H$ such 
 that Kirchhoff's law is satisfied at each vertex of $G$, that is
\[
   \sum_{\vec{e}=u\rightarrow
 v}\phi(\vec{e})+\sum_{\vec{e}=u\leftarrow v}\phi(\vec{e}) =0,
\]
 for every vertex $v$, and where the first sum is taken over all arcs
 towards $v$ and the second sum is over all arcs leaving $v$.  An
 $H$-flow is \emph{nowhere zero} \index{flow!nowhere zero}if $\phi$
 never takes the value 0.

 By replacing the group element on an edge $e$ by its inverse, it is
 clear that two orientations that differ only in the direction of
 exactly one arc $\vec{e}$ have the same number of nowhere zero
 $H$-flows for any $H$. Thus, this number does not depend on the
 choice of orientation of $G$. In fact, when $H$ is finite, it does
 not depend on the 
 structure of the group, but rather only on its cardinality. The
 following, due to Tutte~\cite{Tut54}, relates the number of nowhere
 zero  flows of $G$ over a finite group  and Tutte
 polynomial of $G$.

\begin{theorem}\label{thm:flow_polynomial}
  Let $G=(V,E)$ be a graph and $H$ a finite Abelian group. If
  $\chi^*(G;H)$ denotes the number of nowhere zero $H$-flows then
\[
    \chi^*(G;H) = (-1)^{|E|-r(E)}T(G; 0,1-|H|).
 \]
\end{theorem}
\begin{proof}[sketch]
  Here we use the universality of the Tutte polynomial. If $e$ is an
  ordinary edge of $G$, then the number of nowhere zero $H$-flows in
  $G/e$ can be partitioned into two sets $P_1$ and $P_2$. We let $P_1$
  consist of those that are also nowhere zero $H$-flows in $G\setminus
  e$, and $P_2$ be the complement of $P_1$. Clearly then $|P_1|=
  \chi^*(G\setminus e; H)$. Furthermore, there is a bijection between
  the elements in $P_2$ and the nowhere zero $H$-flows in $G$, and
  thus $|P_2|= \chi^*(G;H)$. It follows that
 \[
  \chi^*(G;H)=\chi^*(G\setminus e;H)-\chi^*(G/e;H),
\]
  and hence $\chi^*(G;H)$ satisfies
  (\ref{eq:TG_invariant_del_con}). It is also easy to check that
  $\chi^*(G;H)$ satisfies (\ref{eq:TG_invariant_prod}). Since
  $\chi^*(L;H)=0$ and $\chi^*(B;H)=|H|-1$, the result follows from
  Theorem~\ref{universal}.  \qed\end{proof} 

 Consequently, $\chi^*(G;\lambda)$ is a polynomial called the
  \emph{flow polynomial} which  for $\lambda$ an integer at least 1
  gives the number of nowhere zero flows of $G$ in a
  group of order $\lambda$. We call any nowhere zero $H$-flow simply a
  $\lambda$-flow if $|H|=\lambda$.

 If the Abelian group is $\mathbb{Z}_3$, and the graph is 4-regular, then the
 Tutte polynomial at $(0,-2)$ counts the number of nowhere zero
 $\mathbb{Z}_3$-flows on $G$. But 
 these flows are in one-to-one correspondence with orientations such
 that at each vertex exactly two edges are directed in and two
 out. Such an orientation is called an \emph{ice configuration}
 \index{configuration!ice} of $G$ (see Lieb~\cite{Lie67} and
 Pauling~\cite{Pau35} for this important model of ice and its physical
 properties). Thus, we have the following corollary.

\begin{corollary}
  If $G$ is a 4-regular graph, then $T(G;0,-2)$ equals the number
  of ice configurations of $G$.
\end{corollary}

 We mentioned that proper colorings are the dual concept of nowhere
  zero flows, and now with Theorem~\ref{flow_polynomial}
  and~(\ref{eq:6::duality}) we observe that 
 \[ \chi(G;\lambda)=\lambda\chi^*(G^*;\lambda),\]
 for $G$ a connected planar graph and $G^*$, any of its plane
 duals. Thus, to each $\lambda$-proper coloring in $G$ corresponds
 $\lambda$ nowhere zero $\mathbb{Z}_{\lambda}$-flows of $G^*$. A
 bijective proof 
 can be found in Diestel~\cite{Die00}.

 Thus, by the 4-color theorem and the duality relation between
 colorings and nowhere zero $H$-flows, every bridgeless planar graph has a
 $4$-flow. For cubic graph, having a nowhere zero
 $\mathbb{Z}_2\times\mathbb{Z}_2$-flow is equivalent to be
 3-edge-colorable. Therefore, as the Petersen graph is not
 3-edge-colorable, it has no  
 $4$-flow. However,  the Petersen graph does have  a
 5-flow. In fact, the famous 5-flow conjecture  of Tutte~\cite{Tut54}
 postulates that every bridgeless graph has a $5$-flow. 

 The 5-flow conjecture is clearly difficult as it is not even apparent
 that every graph will have a $\lambda$-flow for some
 $\lambda$. However, Jaeger~\cite{Jae76} proved that every bridgeless
 graph has an $\mathbb{Z}_2\times\mathbb{Z}_2\times\mathbb{Z}_2$-flow,
 thus every bridgeless graph has a $8$-flow.
 Subsequently Seymour~\cite{Sey81} proved 
 that every bridgeless graph has a
 $\mathbb{Z}_2\times\mathbb{Z}_3$-flow, thus every graph has a $6$-flow.

 Not much is currently known about properties of the flow polynomial
 apart from those that can be deduced from its duality with the
 chromatic polynomial and efforts to solve the 5-flow conjecture.
 However, for some recent work in this direction,
 see Dong and Koh~\cite{DK07} and Jackson~\cite{Jac07}.

\subsection{Abelian Sandpile Models}
\label{subsec:sandpile}

 Self-organized criticality is a concept widely considered in various
   domains since Bak, Tang and Wiesenfeld~\cite{BTW88} introduced
   it. One of the paradigms in this framework is the Abelian sandpile
   model, introduced by Dhar~\cite{Dha90}.

 We begin by recalling the definition of the general Abelian sandpile
   model on a set of $N$ sites labeled 1, 2, $\ldots$, $N$, that we
   refer to as the system. A sandpile at each site $i$ has height 
   given by an integer $h_i$. The set $\vec{h}=\{h_i\}$ is called a
   \emph{configuration}\index{configuration!of a system} of the
   system. For every site $i$, a threshold $H_i$ is defined;
   configurations with $h_i < H_i$ for all $i$ are called
   \emph{stable}\index{configuration!stable}. For every stable
   configuration, the height $h_i$ increases in time at a constant
   rate; this is called the \emph{ loading} of the system. This
   loading continues until $h_i \geq H_i$ for some $i$. The site $i$
   then `topples' and all the values $h_j$, for $1\leq j \leq N$, are
   updated according to the rule:
   \begin{equation}
     \label{eq:toppling}
     h_j= h_j -M_{ij},\ \ \mbox{for all}\ j,
   \end{equation}
  where $M$ is a given fixed integer matrix satisfying
 \[
    M_{ii} > 0,\ \ M_{ij} \leq 0\ \ \mbox{and}\ \
    s_i=\sum_{j}{M_{ij}}\geq 0.
 \] 
  If, after this redistribution, the height at some vertex exceeds its
  threshold, we again apply the toppling rule (\ref{eq:toppling}), and
  so on, until we arrive at a stable configuration and the loading
  resumes. The sequence of topplings is called an
  \emph{avalanche}\index{avalanche}. We assume that an avalanche is
  ``instantaneous'', so that no loading occurs during an avalanche.

  The value $s_i$ is called the \emph{dissipation}\index{dissipation}
  at site $i$. We say that $s_i$ is \emph{dissipative} if $s_i > 0$ and
  \emph{non-dissipative} if $s_i=0$.  It may happen that an avalanche
  continues without end. We can avoid this possibility by requiring
  that from every non-dissipative site $i$, there
  exists a path to a dissipative site $j$.  In
  other words, there is a sequence $i_0,\ldots, i_{n}$, with $i_0=i$,
  $i_n=j$ and $M_{i_{k-1},i_k}<0$, for $k=1,\ldots, n$.  In this
  case, following Gabrielov~\cite{Gab93}, we say that the system is
  \emph{weakly-dissipative}, and we assume that a system is always
  weakly-dissipative.  In a weakly-dissipative system, any
  configuration $\vec{h}$ will eventually arrive at a stable
  configuration. But the process is infinite, and the stable
  configurations are clearly finite. Thus, some stable configurations
  recur, and these are called \emph{critical
  configurations}\index{configuration!critical}.

  The sandpile process has an Abelian property, in that if at some
  stage, two sites can topple, the resulting stable configurations
  after the avalanche is independent of the order in which the sites
  toppled. Thus, for any configuration $\vec{h}$, the process
  eventually arrives at a unique critical configuration $\vec{c}$.

  Let $G$ be a graph, $q\in V(G)$, and $L'$ be the
  minor of the 
  Laplacian of $G$ resulting from deleting the row and column
  corresponding to $q$. When the matrix $M$ is $L'$ for some
  vertex $q$, the Abelian sandpile model coincides with the
  chip-firing game or dollar game on a graph that was defined by
  Biggs~\cite{Big96b}. For the rest of this Subsection we assume
  $M$ is given in this way.

  For a configuration $\vec{h}$, we define its
  \emph{weight}\index{configuration!weight} to be ${w}(\vec{h}) =
  \sum_{i=0}^{N} {h_{i}}$.
  If $\vec{c}$ is a critical configuration, we define its
  \emph{level}\index{configuration!level} as
 \[
      {\rm level}(\vec{c}) =
      {w}(\vec{c})-|E(G)| 
                               +\operatorname{deg}(q). 
 \]
  This definition may seem a little unnatural, but it is justified by
  the following theorem of Biggs~\cite{Big96b}, which tell us that it
  is actually the right quantity to consider if we want to grade the
  critical configurations.
 
  \begin{theorem}\label{2:level}
   If $G=(V,E)$ is a graph and $\vec{h}$ a critical configuration of $G$,
   then
   \[
      0 \leq {\rm level}(\vec{h}) \leq |E|-|V|+1.
   \]
 \end{theorem}
 The right-most quantity is called
  the \emph{cyclomatic number} \index{cyclomatic number} of $G$.   We now consider
  the generating function of the of these critical configurations.

\begin{definition}\index{critical configuration polynomial}\index{polynomial!critical configuration} 
Let $G=(V,E)$ be a graph and for nonnegative integers $i$ let $c_{i}$ be the
  number of critical configurations with level $i$. Then the \emph{critical configuration polynomial} is

  \[
       {P}_{q}(G;y)=
              \sum_{i=0}^{|E|-|V|+1}c_{i}y^{i}. 
  \]

\end{definition}

\begin{theorem}\label{2:main}
  For a graph $G$ and any vertex $q$, the generating function of the
 critical configurations equals the Tutte polynomial of $G$ along the
 line $x=1$, that is,
 \[ 
   {P}_{q}(G; y)= T(G; 1,y),
 \]
 and thus ${P}_{q}(G; y)$ is independent of the choice of
 $q$.
\end{theorem}
 
A proof using deletion and contraction of an edge incident with the
special vertex $q$ can be found in~\cite{Mer97}.
 
New combinatorial identities frequently arise when a new generating
function can be shown to be related to the Tutte polynomial, as in the
following corollary.

\begin{corollary}
If $G$ is a connected graph, then the number of critical
configurations of $G$ is equal to the number of spanning trees, and
the number of critical configurations with level 0 is equal to the
number of acyclic orientations with a unique source.
\end{corollary}

This follows from comparing Theorem \ref{2:main} with Theorems
\ref{trivial_interpretations} and \ref{orientation_results}.

\subsection{The Reliability Polynomial}
\label{subsec:reliability}

  Many of the invariants reviewed thus far have applications in the
  sciences, engineering and computer science, but the
  reliability polynomial we discuss next is among the most directly
  applicable.  
\begin{definition} \index{reliability polynomial} \index{polynomial!reliability}
Let $G$ be a connected graph or network with $n$ vertices and $m$ edges, and suppose that
  each edge is independently chosen to be active with probability
  $p$. Then the (all terminal) \emph{reliability polynomial} is
  
  \begin{equation}
    \begin{split}
      R(G;p)&=\sum_{\substack{A \text{ spanning}\\
      \text{connected}}} p^{|A|}(1-p)^{|E-A|}\\ &=\sum_{k=0}^{m-n+1}
      g_{k}\,p^{k+n-1}(1-p)^{m-k-n+1},
    \end{split}
  \end{equation}
  where $g_k$ is the number of spanning connected subgraphs with
  $k+n-1$ edges.
\end{definition}

Thus the reliability polynomial,
  $R(G;p)$, is the probability that in this random model there is a
  path of active edges between each pair of vertices of $G$.

\begin{theorem}\label{reliability Tutte} 
If G is a connected graph with m edges and n vertices, then
\[
R(G;p)=p^{n-1}(1-p)^{m-n+1}T\left(G;1,\frac{1}{1-p}\right).
\]
\end{theorem}

\begin{proof}  We first note from the rank generating expansion of
  Definition~\ref{def:rank_generating_expansion} that 
\[
T(G;1,y+1)=\sum_{k=0}^{m-n+1} g_{k} y^{k},
\] 
 since the only non-vanishing terms are those corresponding to
 $A\subseteq E$ with $r(E)=r(A)$, that is spanning connected
 subgraphs.

We then observe that
\begin{equation*}
  \begin{split}
    R(G;p)&=\sum_{k=0}^{m-n+1} g_{k} p^{k+n-1}(1-p)^{m-k-n+1}\\
          &=p^{n-1}(1-p)^{m-n+1} \sum_{k=0}^{m-n+1} g_{k}
          \left(\frac{p}{1-p}\right)^{k}\\
          &=p^{n-1}(1-p)^{m-n+1}T\left(G;1,1+\frac{p}{1-p}\right)\\
          &=p^{n-1}(1-p)^{m-n+1}T\left(G;1,\frac{1}{1-p}\right).\\
   \end{split}
\end{equation*}
\qed\end{proof}

If we extend the reliability polynomial to graphs with more than one
component by defining $R(G \cup H;p)=R(G,p)R(H,p)$, then this result
may also be proved using the universality property of the Tutte
polynomial.  Observe that if an ordinary edge is not active (this
happens with probability $1-p$), then the reliability of the network
is the same as if the edge were deleted.  Similarly, if an edge is
active (which happens with probability $p$), then the reliability is
the same as it would be if the edge were contracted.  Thus, the
reliability polynomial has the following deletion/contraction
reduction:

\[
R(G;p) = (1-p)R(G \setminus e)+ pR(G/e).
\]

With this, and noting that $R(G*H;p)=R(G;p)R(H;p)$ with $R(L,p)=1$ and
$R(B;p)=p$,  Theorem~\ref{reliability Tutte} also follows immediately from
Theorem~\ref{universal}.

There is a vast literature about reliability and the reliability
polynomial; for a good survey, including a wealth of open problems, we
refer the reader to Chari and Colbourn~\cite{CC97}.

\subsection{The Shelling Polynomial}
\label{subsec:shelling}

  A {\it simplicial complex}\index{complex!simplicial} $\Delta$ is a
  collection of subsets of a set of vertices $V$ such that if $v \in
  V$, then $\{ v\} \in \Delta$ and also if $F \in \Delta$ and $H
  \subseteq F$, then $H\in \Delta$.  The elements of $\Delta$ are
  called \emph{faces}\index{complex!face}. Maximal faces are called
  \emph{facets}\index{complex!facet}, and if all the facets have the
  same cardinality, $\Delta$ is called
  \emph{pure}\index{complex!pure}. The dimension of a face is its size
  minus one and the dimension of a pure simplicial complex is the
  dimension of any of its facets.

 If $f_{k}$ is the number of faces of size $k$ in a simplicial complex
 $\Delta$, then the vector ($f_0$, $f_1$, \ldots, $f_d$) is called the
 \emph{face vector} or $f$-\emph{vector} \index{f-vector@$f$-vector}
 of $\Delta$, and
\begin{equation}\label{eq:FaceVector}
     {f}_{\Delta}(x) = \sum_{k=0}^{d}{f_k x^{d-k}},
 \end{equation}
 is the generating function of the faces of $\Delta$, or \emph{face
 enumerator}.\index{face enumerator}

 The collection of spanning forests of a connected graph $G$ forms a
 pure $(d-1)$ dimensional simplicial complex $\Delta(G)$.  The points of
 $\Delta(G)$ are the non-loop edges of $G$ and its facets are the spanning
 trees, so $d=r(E)$. The collection of complements of spanning connected
 subgraphs of $G$ also forms a pure $(d^{*}-1)$ dimensional simplicial
 complex 
 $\Delta^{*}(G)$. Here the elements are the non-bridge edges, while the
 facets are complements of spanning trees, when viewed as subsets of
 $E$; in general, if $A$ is the edge-set of a spanning connected
 subgraph of $G$ of cardinality $k+n-1$, then $E\setminus A$ is a face
 of size $m-n+1-k$ in $\Delta^{*}(G)$. Thus, $d^{*}=m-n+1$ and if, as
 before, $g_k$ is the number of spanning connected subgraphs with
 $k+n-1$ edges, the $f$-vector of $\Delta^{*}(G)$ is $(f^*_0,\ldots,
 f^*_{d^{*}})$, where $f^*_i=g_{d^{*}-i}$.

\begin{theorem}\label{TutteFaceVector}

 The Tutte polynomial gives the face enumerators for both $\Delta(G)$
 and $\Delta^{*}(G)$:
 \[
    T(G;x+1,1)=\sum_{k=0}^{d} f_{k} x^{d-k}
                   ={f}_{\Delta(G)}(x),
 \]
 
and
 
\[ 
     T(G;1,y+1)=\sum_{i=0}^{d^{*}} f^*_{i} y^{d^{*}-i}
                    ={f}_{\Delta^*(G)}(x).
 \]
\end{theorem}

\begin{proof} 
This follows readily by comparing (\ref{eq:FaceVector}) with
Definition~\ref{def:rank_generating_expansion}.  \qed\end{proof}

   For a pure simplicial complex $\Delta$, a \emph{shelling}
   \index{shelling} is a linear order of the facets $F_1$, $F_2$,
   $\dots$, $F_t$ such that, if $1 \leq k \leq t$, then $F_k$ meets
   the complex generated by its predecessors, denoted $\Delta_{k-1}$,
   in a non-empty union of maximal proper faces.  A complex is said to
   be \emph{shellable} if it is pure and admits a shelling.  A good
   exposition of the following results can be found in
   Bj\"orner~\cite{Bjo92}.

 For $1 \leq k \leq t$, define $\mathcal{R}(F_k)= \{ x \in F_k |\
   F_k\setminus x \in \Delta_{k-1}\}$,
   where here $\Delta_0=\emptyset$. The number of facets such that
   $|F_k - \mathcal{R}(F_k)|=i$ is denoted by $h_i$ and it does not
   depend on the particular shelling (this follows for example from
   (\ref{eq:face_shelling}) below). The vector ($h_0$, $h_1$,
   $\ldots$, $h_d$) is called the
   $h$--\emph{vector}\index{h-vector@$h$-vector} of $\Delta$.  The
   \emph{shelling polynomial}\index{shelling polynomial} is the
   generating function of the $h$--vector, and is given by
   \[
   {h}_{\Delta}(x) = \sum_{i=0}^{d}{h_i x^{d-i}}.
   \] 
   
  The face enumerator and shelling polynomial are related in a
   somewhat surprising way, namely
   \begin{equation}\label{eq:face_shelling}
     {h}_{\Delta}(x+1) = {f}_{\Delta}(x).
   \end{equation} 
   
   Both $\Delta(G)$ and $\Delta^{*}(G)$ are known to be shellable,
   see for example  Provan and Billera~\cite{PB80}), and
   thus~(\ref{eq:face_shelling})  
   gives the following corollary to Theorem~\ref{TutteFaceVector},
   relating the two shelling polynomials to the Tutte polynomial
   (see Bj\"orner~\cite{Bjo92}).
  
\begin{corollary}\label{shelling Tutte} Let $G$ be a graph.  Then
   \[
   T(G; x,1)= h_{\Delta(M)}(x) = \sum_{i=0}^{d}{h_i x^{d-i}}
   \]
   and
   \[
   T(G; 1,y)= h_{\Delta*(G)}(y)= \sum_{i=0}^{d^{*}}{h_i^{*}
   y^{d^{*}-i}}.
   \]
 \end{corollary}
  
   The reader may have noticed that the reliability polynomial as well
   as the face enumerator and shelling polynomial of $\Delta^{*}(G)$
   are all specializations of the Tutte polynomial along the line
   $x=1$. There is an important open conjecture in algebraic
   combinatorics about the $h$-vectors (and hence the shelling
   polynomials), of the two complexes coming from a graph (or, more
   generally, a matroid), namely that they are `pure O-sequences'. For
   more details see Stanley~\cite{Sta96b} or \cite{Mer01}.  The latter
   also relates the shelling polynomial and the chip firing game. Let
   $G$ be a graph with $n$ vertices and $m$ edges. From
   Corollary~\ref{shelling Tutte} and Theorem~\ref{2:main}, we get
   that $c_i=h_{m-n+1-i}^{*}$, where $c_i$ is the number of critical
   configurations of level $i$ of $G$ and ($h_{0}^{*}$, \ldots,
   $h_{m-n+1}^{*}$) is the $h$-vector of
   $\Delta^{*}(G)$. In~\cite{Mer01} it is proved that ($c_{m-n+1}$,
   \ldots, $c_0$) is a pure O-sequence. Thus, the conjecture is true
   for the simplicial complex $\Delta^{*}(G)$ but is still open for
   $\Delta(G)$.  

It is also clear from Theorem~\ref{reliability Tutte} and
Corollary~\ref{shelling Tutte} that the reliability and shelling
polynomials are related. This connection is explored, and open
questions related to it presented, by Chari and Colbourn~\cite{CC97}.

\section{Some Properties of the Tutte Polynomial}
\label{sec:properties}

There is a large and ever-growing body of information about properties
of the Tutte polynomial.  Here, we present some of them, again with an
emphasis on illustrating general techniques for extracting information
from a graph polynomial.

\subsection{The Beta Invariant}

  Even a single coefficient of a graph polynomial can encode a
  remarkable amount of information.  It may characterize entire
  classes of graphs and have a number of combinatorial
  interpretations.  A noteworthy example is the $\beta$
  invariant\index{beta invariant@$\beta$ invariant}, introduced (in
  the context of matroids) by Crapo in~\cite{Cra67}.

  \begin{definition}
   Let $G=(V,E)$ be a graph with at least two edges.  The $\beta$ invariant
   of $G$ is
   \[
     \beta \left( G \right) = \left( { - 1} \right)^{r\left( G
   \right)} \sum\limits_{A \subseteq E} {\left( { - 1} \right)^{\left|
   A \right|} r\left( A \right)}.
    \]   
   \end{definition}

  The beta invariant is a deletion/contraction invariant, that is, it
  satisfies~(\ref{eq:TG_invariant_del_con}). However, the $\beta$
  invariant is zero if and only if $G$ either has loops or is not
  two-connected.  Thus, the $\beta$ invariant  is  not a
  Tutte-Gr\"othendieck invariant in the sense of Section~\ref{sec:IV}.
  While the $\beta$ invariant may be 
  defined to be 1 for a single edge or a single loop, it  still 
   will not satisfy~(\ref{eq:TG_invariant_prod}), and it is not
  multiplicative with respect to disjoint unions and one-point joins.
  Nevertheless, the $\beta$ invariant derives from the Tutte
  polynomial.

  \begin{theorem}
   If $G$ has at least two edges, and we write $T(G;x,y)$ in the form
   $\sum{t_{ij} x^i y^j}$, then $t_{0,1} = t_{1,0} $, and this common
   value is equal to the $\beta$ invariant.
  \end{theorem}

   \begin{proof}  This can easily be proved by induction, using
     deletion/contraction for an ordinary edge, and otherwise noting
     that the $\beta$ invariant is zero if the graph has loops or is
     not two-connected.  \qed\end{proof}

   The $\beta$ invariant does not change with the insertion of
   parallel edges or edges in series.  Thus, homeomorphic graphs have
   the same $\beta$ invariant. The $\beta$ invariant is also
   occasionally called the chromatic invariant, because $\chi '\left(
   {G;1} \right) = \left( { - 1} \right)^{r\left( G \right)} \beta
   \left( G \right)$, where $\chi(G;x)$ is the chromatic polynomial.
 
 \begin{definition}
   A series-parallel\index{graph!series-parallel} graph is a graph
   constructed from a digon (two vertices joined by two edges in
   parallel) by repeatedly adding an edge in parallel to an existing
   edge, or adding an edge in series with an existing edge by
   subdividing the edge. Series-parallel graphs are
   loopless multigraphs, and are planar.
 \end{definition} 

  Brylawski, \cite{Bry71} and also~\cite{Bry82}, in the context of
  matroids, showed that the $\beta$ invariant completely characterizes
  series-parallel graphs.

  \begin{theorem}
   $G$ is a series-parallel graph if and only if $\beta \left( G
  \right) = 1$.
  \end{theorem} 

  Using the deletion/contraction definition of the Tutte polynomial,
  it is quite easy to show that the $\beta$ invariant is unchanged by
  adding an edge in series or in parallel to another edge in the
  graph.  This, combined with the $\beta$ invariant of a digon being
  one, suffices for one direction of the proof. The difficulty is in
  the reverse direction, and the proof is provided in~\cite{Bry72} by
  a set of equivalent characterizations for series-parallel graphs,
  one by excluded minors and another that the $\beta$ invariant is 1
  for series-parallel graphs.  For graphs, the excluded minor is $K_4$
  (cf. Duffin~\cite{Duf65} and Oxley~\cite{Oxl82}).  Succinct proofs may
  also be found in Zaslavsky~\cite{Zas87}. The fundamental
  observation, which may be applied to other situations, is that there
  is a graphical element, here an edge which is in series or parallel
  with another edge, which behaves in a tractable way with respect to
  the computation methods of the polynomial.

  The $\beta$ invariant has been explored further, for example by
  Oxley in~\cite{Oxl82} and by Benashski, Martin, Moore and Traldi
  in~\cite{BMMT95}.  Oxley characterized 3-connected matroids with
  $\beta \leqslant 4$, and a complete list of all simple 3-connected
  graphs with $\beta \leqslant 9$ is given in~\cite{BMMT95}.

  A wide variety of combinatorial interpretations have also been found
  for the $\beta$ invariant. Most interpretations involve objects
  other than graphs, but we give two graphical interpretations
  below. The first is due to Las Vergnas~\cite{Las84}.

  \begin{theorem}
    Let $G$ be a connected graph.  Then $2\beta \left( G \right)$
    gives the number of orientations of $G$ that have a unique source
    and sink, independent of their relative locations.
  \end{theorem}

  This result is actually a consequence of a more general theorem
  giving an alternative formulation of the Tutte polynomial, which
  will be discussed further in Subsection~\ref{subsec:coefficient_relations}. We also have the
  following result from~\cite{E-M04}.
 
 \begin{theorem}
    Let $G=(V,E)$ be a connected planar graph with at least two edges.
    Then
    \[
       \beta= \frac{1}{2}\sum {\left( { - 1} \right)^{c\left(
          E\setminus P \right) + 1} },
     \]
    where the sum is over all closed trails $P$ in $\vec G_m $ which
    visit all its vertices at least once.
  \end{theorem}

 Like the interpretations for $T(G;x,x)$ given in the
 Subsection~\ref{subsec:x=y}, this result follows from the Tutte
 polynomial's relation to the Martin polynomial.

  Graphs in a given class may have $\beta$ invariants of a particular
  form.  McKee~\cite{McK01} provides an example of this in
  dual-chordal graphs.  A \emph{dual-chordal
  graph}\index{graph!dual-chordal} is 2-connected, 3-edge-connected,
  such that every cut of size at least four creates a bridge. A
  \emph{$\theta$ graph} has two vertices with three edges in parallel between
  them.  A dual-chordal graph has the property that it may be reduced
  to a $\theta$ graph by repeatedly contracting induced subgraphs of
  the following forms: digons, triangles, and $K_{2,3}$'s, where in
  all cases each vertex has degree 3 in $G$.
  
  \begin{theorem}
  If $G$ is a dual-chordal graph, then $\beta(G)=2^a5^b$. Here, $a$ is
  the number of triangles in $G$, where each vertex has degree 3 in
  $G$, that are contracted in reducing $G$ to a $\theta$ graph.
  Similarly, $b$ is the number of induced $K_{2,3}$ in $G$, again
  where each vertex has degree 3 in $G$, that are contracted in
  reducing $G$ to a $\theta$ graph.
  \end{theorem}
   
 The proof follows from considering the acyclic orientations of $G$
 with unique source and sink and applying the results of Green and
 Zaslavsky~\cite{GZ83}.

 \subsection{Coefficient Relations}
\label{subsec:coefficient_relations}

  After observing that $t_{1,0} = t_{0,1} $ in the development of the
  $\beta$ invariant, it is natural to ask if there are similar
  relations among the coefficients $t_{ij} $ of the Tutte polynomial
  $T(G;x,y) = \sum {t_{ij} x^i y^j } $ and whether there are
  combinatorial interpretations for these coefficients as well.  The
  answer is yes, although less is known.  The most basic fact, and one
  which is not obvious from the rank-nullity formulation of
  Definition~\ref{def:rank_generating_expansion}, is that all the
  coefficients of the Tutte polynomial are non-negative.

  That $t_{1,0} = t_{0,1} $ is one of an infinite family of relations
  among the coefficients of the Tutte polynomial.
  Brylawski~\cite{Bry82} has shown the following:

  \begin{theorem}\label{coeff_rel}
    If $G$ is a graph with at least $m$ edges, then
    \[
    \sum_{i=0}^{k} {\sum_{j=0}^{k-i} 
          {(-1)^j
          \binom{k-i}{j}  t_{ij} }}= 0,
   \] 
  for $k=0,1\ldots, m-1$.
  \end{theorem}
  

  Additionally, Las Vergnas in~\cite{Las84}  found combinatorial
  interpretations in the context of oriented matroids 
  for these coefficients by determining yet another generating function
  formulation for the Tutte polynomial.  Gioan and Las
  Vergnas~\cite{GLV05} give the following specialization to orientations of
  graphs.  

  \begin{theorem}
     Let $G$ be a graph with a linear ordering of its edges. Let $o_{i,j}$
     be the number of orientations of G such that the number of edges that
     are smallest on some consistently directed cocycle is i and the number
     of edges that are smallest on a consistently directed cycle is j.  Then
     \[
      T(G;x,y) = \sum\limits_{i,j} {o_{i,j} 2^{-(i+j)} x^i y^j }
      ,
    \]
     and thus $t_{ij} = o_{i,j}/(2^{i + j})$.
  \end{theorem}

  The proof is modeled on Tutte's proof that the $t_{ij}$'s are
  independent of the ordering of the edges by using
  deletion/contraction on the greatest edge in the ordering. 

\begin{figure}[hbtp]
      \begin{center}
    \includegraphics[scale=0.75]{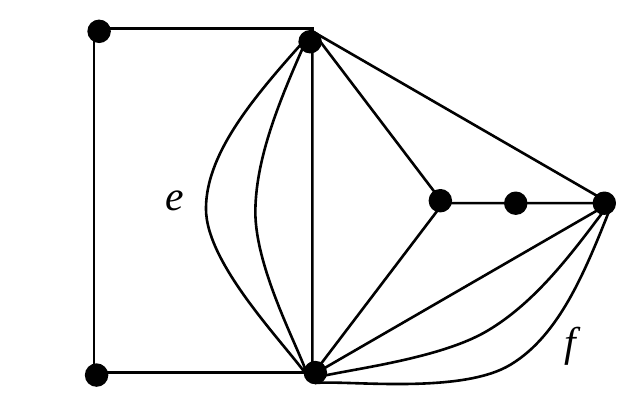}
      \caption{A counterexample to the conjecture}\label{Fig:counterexample} 
      \end{center}
      \end{figure}

  Another natural question is to ask if these coefficients are
  unimodular or perhaps log concave, for example in either $x$ or $y$.
  While this was originally conjectured to be true (see Seymour and
  Welsh~\cite{SW75}, Tutte~\cite{Tut84}), then Schw\"arzler~\cite{Sch93}
  found a contradiction in the graph in
  Fig.~\ref{Fig:counterexample}. This counterexample can be extended
  to an infinite family of counterexamples by increasing the number of
  edges parallel to $e$ or $f$.

  The unimodularity question for the chromatic polynomial, raised
  by Read in~\cite{Rea68}, is still unresolved. 

 \subsection{Zeros of the Tutte Polynomial}
\label{subsec:zeros}
  Because the Tutte polynomial is after all a polynomial, it is very
  natural to ask about its zeros and factorizations.  The importance of its zeros is
  magnified by their interpretations.  For example, since $T(G;0,y)$
  is essentially the flow polynomial, a root of the form $(0,
  1-\lambda)$, for $\lambda$ a positive integer, means that $G$ does
  not have a nowhere zero flow for any Abelian group of order
  $\lambda$.  Similarly, since $T(G;x,0)$ is essentially the
  chromatic polynomial, a root of the form $( 1-\lambda,0)$ with
  $\lambda$ a positive integer, means that $G$ cannot be properly
  colored with $\lambda$ colors.  In particular, a direct proof the
  four-color theorem would follow if it could be shown that the Tutte
  polynomial has no zero of the form $(-3,0)$ on the class of planar
  graphs.  Of course, because of the duality between the flow and
  chromatic polynomials, results for the zeros of the one informs the
  other, and vice versa.  Jackson~\cite{Jac03} surveys zeros of both
  chromatic and the flow polynomials.
  
  As we will see in detailed in the next chapter, the chromatic
  polynomial has an additional interpretation as the zero-temperature
  antiferromagnetic Potts model of statistical mechanics.  In this
  context, its zeros correspond to  numbers of spins for which the ground
  state degeneracy function may be nonanalytic. This has led to
  research into its zeros by theoretical physicists as well as
  mathematicians.  Traditionally, the focus from a graph theory
  perspective was on positive integer roots of the chromatic
  polynomial, corresponding to graph not being properly colorable with
  $q$ colors.  In statistical mechanics however, the relevant quantity
  involves the limit of an increasing family of graphs as the number,
  $n$, of vertices goes to infinity.  This shifted the focus to the
  complex roots of the chromatic polynomial, since the sequence of
  complex roots as $n \to \infty $ may have an accumulation point on
  the real axis.
  
  Because of this, a significant body of work has emerged in recent
  years devoted to clearing regions of the complex plane (in
  particular regions containing intervals of the real axis) of roots
  of the chromatic polynomial.  Results showing that certain intervals
  of the real axis and certain complex regions are free of zeros of
  chromatic polynomials include those of Woodall~\cite{Woo92},
  Jackson~\cite{Jac93}, Shrock and Tsai~\cite{ST97a,ST97b},
  Thomassen~\cite{Tho97}, Sokal \cite{Sok01b}, Procacci, Scoppola, and
  Gerasimov~\cite{PSG03}, Choe, Oxley, Sokal, and
  Wagner~\cite{COSW04}, Borgs~\cite{Bor06}, and Fernandez and
  Procacci~\cite{FP}.  One particular question concerns the maximum
  magnitude of a zero of a chromatic polynomial and of zeros
  comprising region boundaries in the complex plane as the number of
  vertices $n \to \infty $.  An upper bound is given in~\cite{Sok01b},
  depending on the maximal vertex degree.  There are, however,
  families of graphs where both of these magnitudes are unbounded (see
  Read and Royle~\cite{RR91}, Shrock and Tsai~\cite{ST97a,ST98},
  Brown, Hickman, Sokal and Wagner~\cite{BHSW01}, and
  Sokal~\cite{Sok04}).  For recent discussions of some relevant
  research directions concerning zeros of chromatic polynomials and
  properties of their accumulation sets in the complex plane, as
  well as approximation methods, see,  e.g.,  Shrock and
  Tsai~\cite{ST97b}, Shrock~\cite{Shr01}, 
  Sokal~\cite{Sok01a, Sok01b}, Chang and Shrock~\cite{CS01b}, Chang,
  Jacobsen, Salas, and Shrock~\cite{CJSS04}, Choe, Oxley, Sokal, and
  Wagner~\cite{COSW04}, Dong and Koh~\cite{DK04}, and more recently
  Royle~\cite{Roya, Royb}.
 
 If $G$ is  a graph with chromatic number $k+1$, then  $\chi(G;x)$ has
 integer roots at $0,1,\ldots, 
 k$. Thus, the chromatic polynomial of $G$ can be written as
\[ 
    \chi(G;x)=x^{a_0}(x-1)^{a_1}\cdots(x-k)^{a_{k}}q(x),
\]
 where $a_0,\ldots,a_{k}$ are integers and $q(x)$ is a polynomial
 with no integer roots in the interval $[0,k]$.
 In contrast to this we have the following result of Merino, de Mier and
 Noy~\cite{MMN01}.

\begin{theorem}\label{irreducible}
 If $G$ is a 2-connected graph, then $T(G;x,y)$ is irreducible in
 $\mathbb{Z}[x,y]$. 
\end{theorem}
 The proof is quite technical and it heavily relies on
 Theorem~\ref{coeff_rel} and that $\beta(G)\neq 0$ if and only if $G$
 has no loops and it is 2-connected.

 If $G$ is not 2-connected, then $T(G;x,y)$ can be factored. From
 Proposition~\ref{product-formula} we get that if $G$ is a 
 disconnected graph with connected components $G_1, \ldots, 
 G_{\kappa}$, then $T(G;x,y)=\prod_{i=1}^{\kappa}T(G_i;x,y)$. So
 let us consider when $G$ is connected but not 2-connected.

 One of the basic properties mentioned in~\cite{BO92} is that
  $y^{s}|T(G;x,y)$ if and only if $G$ has $s$ loops. 
 Thus, let us focus on  loopless connected graphs that are not
  2-connected.  It is well-known
 that such graphs have a decomposition
 into its blocks, see for example~\cite{Bol98}. A \emph{block of}
 \index{graph!block of} a
 graph $G$ is either a bridge or a maximal 2-connected subgraph. If
 two blocks of $G$ intersect, they do so in a cut vertex. 
 By Theorem~\ref{irreducible} and Proposition~\ref{product-formula} we
 get the following. 

\begin{corollary}
   If $G$ is a loopless connected graph that is not 2-connected with
    blocks $H_1,\ldots, H_p$, then the factorization of
    $T(G;x,y)$ in  $\mathbb{Z}[x,y]$ is exactly
\[
   T(G;x,y)=T(H_1;x,y)\cdots T(H_p;x,y).
\]     
\end{corollary}

 \subsection{Derivatives of the Tutte Polynomial}
\label{subsec:derivatives}

  It is also most natural to differentiate the Tutte polynomial and to
  ask for combinatorial interpretations of its derivatives.  For
  example, Las Vergnas~\cite{Las} has found the following
  combinatorial interpretation of the derivatives of the Tutte
  polynomial.  It first requires a slight generalization of the
  notions of internal and external activities given in
  Subsection~\ref{subsec:trees_expansion}.
  
  \begin{definition}\index{active edge!externally}
                    \index{active edge!internally} 
  Let $G=(V,E)$ be a graph with a linear order on its edges, and let $A
  \subseteq E$.  An edge $e \in A$ and a cut $C$
  are internally active with respect to $A$ if $e \in C \subseteq
  \left( {E\backslash A} \right) \cup \{ e\} $ and $e$ is the smallest
  element in $C$.  Similarly, an edge $e \in E\backslash A$ and a
  cycle $C$ are externally active with respect to $A$ if $e \in C
  \subseteq A \cup \{ e\} $.
 \end{definition}
  
  In the case that $A$ is a spanning tree, this reduces to the
  previous definitions of internally and externally active.
  
  \begin{theorem}
 Let $G$ be a graph with a linear ordering on its edges.  Then
   \[
     \frac{{\partial ^{p + q} }} {{\partial x^p \partial y^q }}
     T(G;x,y) = p!\,q!\,\sum {x^{\operatorname{in}(A)}
     y^{\operatorname{ex}(A)} },
   \]
  where the sum is over all subsets $A$ of the edge set of $G$ such that
  $r\left( G \right) - r\left( A \right) = p$ and $\left| A \right| -
  r\left( A \right) = q$, and where $\operatorname{in}( A )$ is the
  number of internally active edges with respect to $A$, and
  $\operatorname{ex}(A)$ is the number of externally active edges with
  respect to $A$.
  \end{theorem}
  
  The proof begins by differentiating the spanning tree definition of
  the Tutte polynomial, Definition~\ref{def:trees_expansion}, which
  gives a sum over $i$ and $j$ restricted by $p$ and $q$.  This is
  followed by showing that the coefficients of $x^{i - p} y^{j - q}$
  enumerate the edge sets described in the theorem statement.  The
  enumeration comes from examining, for each subset $A$ of $E$, the
  set of $e \in E\backslash A$ such that there is a cut-set of $G$
  contained in $E\backslash A$ with $e$ as the smallest element (and
  dually for cycles).
  
  The Tutte polynomial along the line $x=y$ is a polynomial in one
  variable that, for planar graphs, is related to the Martin
  polynomial via a medial graph construction. From this relationship,
  \cite{E-M04b} derives an interpretation for the $n$-th derivative of
  this one variable polynomial evaluated at $(2,2)$ in terms of edge
  disjoint closed trails in the oriented medial graph.

   \begin{definition}
  For an oriented graph $\vec{G}$ , let $P_n$ be the set of ordered
        $n$-tuples ${\bar p}:=(p_1 , \ldots ,p_n)$, where the $p_i$'s are
        consistently oriented edge-disjoint closed trails in
        $\vec{G}$. 
   \end{definition}

  
  \begin{theorem}
  If G is a connected planar graph with oriented medial graph
  $\overrightarrow {G_m} $, then, for all non-negative integers n,
 \[
  \left. \frac{{\partial ^{n} }} {{\partial x}}
     T(G;x,x)\right|_{x=2} = \sum_{k = 0}^n {(- 1)^{n - k}
     \frac{{n!}}{{k!}}\sum_{\bar p \in P_k \left( {\vec G_m }
     \right)}^{} {2^{m\left( {\bar p} \right)} } }, \]
 where $m\left( {\bar p} \right)$ is the number of vertices of
 $\vec{G}$  not belonging to any of the trails in $\bar p$.
 \end{theorem}

\subsection{ Convolution and the Tutte Polynomial}
\label{subsec:convolution}

  Since the Tutte polynomial can also be formulated as a generating
  function, the tools of generating functions, such as M\"obius
  inversion and convolution, are available to analyze it. A
  comprehensive treatment of convolution and M\"obius inversion can be
  found in Stanley~\cite{Sta96}. Convolution
  identities are valuable because they write a graph polynomial in
  terms of the polynomials of its substructures, thus facilitating
  induction techniques. We have the following result from Kook,
  Reiner, and Stanton~\cite{KRS99} using this approach.
  
   \begin{theorem} The Tutte polynomial can be expressed as
  \[
   T(G;x,  y) = \sum{T(G/A;x,0)
                    T\left({\left. G \right|_A;0,y}\right)},
  \]
   where the sum is over all subsets $A$ of the edge set of G, and
  where $\left. G \right|_A $ is the restriction of G to the edges of
  $A$, i.e. $\left. G \right|_A=G\setminus (E\setminus A)$.
  \end{theorem}

  This result is particularly interesting in that it essentially
  writes the Tutte polynomial of a graph in terms of the chromatic and
  flow polynomials of its minors.  It may be proved in several ways,
  for example by induction using the deletion/contraction relation, or
  from the spanning trees expansion of the Tutte polynomial.  However,
  we present the first proof from~\cite{KRS99} to illustrate the
  technique, which is dependent on results of Crapo~\cite{Cra69}. 
  
  \begin{proof}[sketch]
  We begin with a convolution product of two functions on graphs into
  the 
  ring $\mathbb{Z}[x,y]$ given by \index{convolution}
  $f * g$ $=$ $\sum_{A \subseteq E\left( G \right)} {f\left( {\left. G
  \right|_A } \right)} g\left( {G/A} \right)$. The identity for
  convolution is $\delta (G)$ which is 1 if and only if $G$ is
  edgeless and 0 otherwise.  From Crapo~\cite{Cra69}, we have that
  \[
  T(G;x + 1,y + 1) = \left( {\zeta (1,y) * \zeta (x,1)}
  \right)\left( G \right),
  \]
  where $\zeta \left( {x,y} \right)\left( G \right) = x^{r\left( G
  \right)} y^{r\left( {G^ * } \right)} $. Kook, Reiner, and
  Stanton~\cite{KRS99} then show that $\zeta \left( {x,y} \right)^{ -
  1} = \zeta \left( { - x, - y} \right)$.  From this it follows that
  $T(G;x + 1,0) = ( \zeta (1, - 1) * \zeta (x,1))(G)$ and
  $T(G;0,y + 1) = (\zeta (1,y) * \zeta ( - 1,1))( G )$.  Thus,
  $\sum T(\left. G \right|_A ;0,y + 1)\, T(G/A;x + 1,0) =
  (\zeta (1,y) * \zeta ( - 1,1)) * (\zeta (1, - 1) * \zeta(x,1))( G
  )$. By associativity, the last expression is the same as $(\zeta
  (1,y) * (\zeta ( - 1,1) * \zeta (1, - 1)) * \zeta (x,1))( G ) =
  (\zeta ( 1,y) * \zeta (x,1))( G )= T(G;x + 1,y + 1)$.
  \qed\end{proof}
  
  A formula, known as Tutte's identity for the chromatic polynomial,
  with a similar flavor, exists for the chromatic polynomial.
  
   \begin{theorem}
  The chromatic polynomial can be expressed as
  \[
   \chi(G;x + y) = \sum {\chi(\left. G \right|_A;x)\, \chi(\left. G
  \right|_{A^c};y)},
  \]
  where the sum is over all subsets $A$ of the set of vertices of $G$,
  and where $\left. G \right|_A $ is the restriction of $G$ to the
  vertices of $A$.
  \end{theorem}

  \begin{proof}
  Consider an $(m + n)$-coloring of $G$, and let $A$ be the vertices
  colored by the first $m$ colors.  Then an $(m + n)$-coloring of $G$
  decomposes into an $m$ coloring of $\left. G \right|_A $ using the
  first $m$ colors and an $n$ coloring of $\left. G \right|_{A^c } $
  using the remaining colors.  Thus, for any two non-negative integers
  $m$ and $n$, it follows that $\chi(G;m + n) = \sum {\chi(\left. G
  \right|_A ;m)\,\chi(\left. G \right|_{A^c } ;n )} $.  Since the
  expressions involve finite polynomials, this establishes the result
  for indeterminates $x$ and $y$. \qed
  \end{proof}

\section{The Complexity of the Tutte Polynomial}
\label{sec:complexity}

We assume the reader is familiar with the basic notions of
computational complexity, but for formal definitions in the present
context, see, for example,
Welsh~\cite{Wel93}.

 We have seen that along different algebraic curves in the $XY$ plane,
the Tutte polynomial evaluates to many diverse quantities. Some of
these, such as $T(G;2,2)=2^{|E|}$ are very easy to compute, and others
such as $T(G;1,1)$ may also be computed efficiently, as in
Subsection~\ref{subsec:trees}. In general though, the Tutte polynomial
is intractable, as shown in the following theorem of Jaeger, Vertigan
and Welsh~\cite{JVW90}.

 \begin{theorem}\label{teo9.1}
 The problem of evaluating the Tutte polynomial of a graph at a point
 $(a,b)$ is $\# P$-hard except when $(a,b)$ is on the special
 hyperbola
 \[
         H_1 \equiv (x-1)(y-1) = 1
 \]
 or when $(a,b)$ is one of the special points $(1,1)$, $(-1,-1)$,
 $(0,-1)$, $(-1,0)$, $(i,-i)$, $(-i,i)$, $(j,j^2)$ and $(j^2,j)$,
 where $j = e^{2\pi i/3}$.  In each of these exceptional cases the
 evaluation can be done in polynomial time.
 \end{theorem}

 For planar graphs there is a significant difference.  The technique
 developed using the Pfaffian to solve the Ising problem for the plane
 square lattice by Kasteleyn \cite{Kas61} can be extended to give a
 polynomial time algorithm for the evaluation of the Tutte polynomial
 of any planar graph along the special hyperbola
\[
     H_2 \equiv (x-1)(y-1) = 2.
\]

However, even restricting a class of graphs to its planar members, or
further restricting colouring enumeration on the square lattice, does
not necessarily yield any additional tractability, as shown by the
following results, the first due to Vertigan and Welsh~\cite{VW92},
and the second to Farr~\cite{Far06}.

 \begin{theorem}\label{teo9.2}
 The evaluation of the Tutte polynomial of bipartite planar graphs at
 a point $(a,b)$ is $\# P$-hard except when
 \[
  (a,b) \in H_1 \cup H_2 \cup \{ (1,1),(-1,-1),(j,j^2),(j^2,j)\},
 \]
at which points it is computable in polynomial time.
 \end{theorem}

\begin{theorem}
For $\lambda\geq 3$, computing the number of $\lambda$-colorings of
induced subgraphs of the square lattice is $\# P$-complete.
\end{theorem}

A natural question then arises as to how well an evaluation of the
Tutte polynomial might be approximated. That is, if there is a fully
polynomial randomized approximation scheme, or FPRAS, for $T$ at
a point $(x,y)$ for a well-defined family of graphs. Here, FPRAS
refers to a probabilistic algorithm that takes the input $s$ and the
degree of accuracy $\epsilon$ to produce, in polynomial time on $|s|$
and $\epsilon^{-1}$, a random variable which approximates
$T(G;x,y)$ within a ratio of $1+\epsilon$ with probability
greater than or equal to 3/4. For example, Jerrum and
Sinclair~\cite{JS93} show that there exits an FPRAS for $T$ along
the positive branch of the hyperbola $H_{2}$.

However, in general approximating is provably difficult as well.
Recently, Goldberg and Jerrum~\cite{GJ07} have extended the region of
the $x$-$y$ plane for which the Tutte polynomial does not have an
FPRAS, to essentially all but the first quadrant (under the assumption
that $RP\neq NP$).  A consequence of this is that there is no FPRAS
for counting nowhere zero $\lambda$-flows for $\lambda > 2$. They also
provide a good overview of prior results. For a somewhat more
optimistic prognosis in the case of dense graphs, we refer the reader
to ~\cite{WM00}, and to Alon, Frieze, and Welsh~\cite{AFW95}.

There has been an increasing body of work since the seminal results of
Robertson and Seymour~\cite{RS83, RS84, RS86} impacting computational
complexity questions for graphs with bounded tree-width (see
Bodlaender's accessible introduction in~\cite{Bod93}).  A powerful
aspect of this work is that many $NP$-Hard problems become tractable
for graphs of bounded tree-width. For example, Noble~\cite{Nob98} has
shown that the Tutte polynomial may be computed in polynomial time (in
fact requires only a linear number of multiplications and additions)
for rational points on graphs with bounded tree width.  Makowsky,
Rotics, Averbouch and Godlin~\cite{M+06} provide similar results for
bounded clique-width (a notion with significant computation complexity
consequences analogous to those for bounded tree-width -- see Oum and
Seymour~\cite{OS06}).  Noble~\cite{Nob07} gives a recent survey of complexity results for this area, including new monadic second order logic methods and extensions to the multivariable generalizations of the Tutte polynomial discussed in the next chapter.

Although the Tutte polynomial is not in general computationally
tractable, there are some resources for reasonably sized graphs (about
100 edges).  These include Sekine, Imai, and Tani~\cite{SIT95}, which
provides an algorithm to implement the recursive definition. Common
computer algebra systems such as Maple and Mathematica will compute
the Tutte polynomial for smallish graphs, 
and there are also some implementations freely available in the Web,
such as \url{http://ada.fciencias.unam.mx/~rconde/tulic/} by R. Conde
or \url{http://homepages.mcs.vuw.ac.nz/~djp/tutte/} by G. Haggard and
D. Pearce. 


\section{Conclusion}
\label{sec:conclusion}

 For further exploration of the Tutte polynomial and its properties,
 we refer the reader to the relevant chapters of Welsh~\cite{Wel93}
 and 
 Bollob\'as~\cite{Bol98} for excellent introductions,  and to
 Brylawki~\cite{Bry82}, Brylawski and Oxley~\cite{BO92}, and
 Welsh~\cite{Wel99} for an in-depth treatment of the Tutte polynomial,
 including generalizations to matroids. Although we focused on graphs
 here to broaden accessibility, matroids, rather than graphs, are the
 natural domain of the Tutte polynomial, and Crapo~\cite{Cra69} gives
 a compelling justification for this viewpoint.  Farr~\cite{Far07}
 gives a recent treatment and engaging history of the Tutte
 polynomial. Finally, we especially recommend Tutte's own account of
 how he ``...became acquainted with the Tutte polynomial...''
 in~\cite{Tut04}.

%
%
%

%
%



\printindex
\end{document}